# Asymptotic properties of false discovery rate controlling procedures under independence

**Pierre Neuvial***

*Laboratoire de Probabilités et Modèles Aléatoires*
*Université Paris Diderot — Paris 7*
*175 rue du Chevaleret, 75013 Paris, France*

*INSERM, U900, Paris, F-75248 France*
*École des Mines de Paris, ParisTech, Fontainebleau, F-77300 France*
*Institut Curie, 26 rue d'Ulm, Paris cedex 05, F-75248 France*
*e-mail:* pierre.neuvial@curie.fr

**Abstract:** We investigate the performance of a family of multiple comparison procedures for strong control of the False Discovery Rate (FDR). The FDR is the expected False Discovery Proportion (FDP), that is, the expected fraction of false rejections among all rejected hypotheses. A number of refinements to the original Benjamini-Hochberg procedure [1] have been proposed, to increase power by estimating the proportion of true null hypotheses, either implicitly, leading to one-stage adaptive procedures [4, 7] or explicitly, leading to two-stage adaptive (or plug-in) procedures [2, 21].

We use a variant of the stochastic process approach proposed by Genovese and Wasserman [11] to study the fluctuations of the FDP achieved with each of these procedures around its expectation, for independent tested hypotheses.

We introduce a framework for the derivation of generic Central Limit Theorems for the FDP of these procedures, characterizing the associated regularity conditions, and comparing the asymptotic power of the various procedures. We interpret recently proposed one-stage adaptive procedures [4, 7] as fixed points in the iteration of well known two-stage adaptive procedures [2, 21].

## Contents

*This work was supported by the association Courir pour la vie, courir pour Curie, and by the ANR project TAMIS.

## 1. Introduction

Multiple testing problems arise when many binary tests are performed simultaneously. The rejection of all hypotheses with individual $p$-values smaller than a fixed threshold results in an increasing number of false discoveries as the number of tests increases. There is therefore a need for a risk measure taking multiple testing into account. Multiple testing procedures use a collection of $p$-values as input and output a set of hypotheses to be rejected.

Since the seminal article by Benjamini and Hochberg [1], the False Discovery Rate (FDR) has been accepted as a practical and convenient method for risk assessment in multiple testing problems involving high-dimensional data analysis, including non parametric estimation by wavelet methods in image analysis, functional magnetic resonance imaging (fMRI) in medicine, source detection in astronomy, or DNA microarray analysis in biology. The FDR is the expected proportion of erroneous rejections among all rejections.

The procedure originally proposed by Benjamini and Hochberg [1], which we term procedure BH95, controls FDR when the true null hypotheses are independent, or display certain forms of positive dependence [1, 3]. Considerable efforts have been devoted to the development of procedures retaining the FDR-controlling capabilities of the BH95 procedure under more general conditions of dependence [3, 19], and/or with higher power [2, 4, 7, 21, 23]. The second aim is driven by the observation that the original BH95 procedure actually controls FDR at level of exactly $\pi_0\alpha$, where $\pi_0$ is the (unknown) proportion of true null hypotheses [3, 8, 19, 23]. When $\pi_0 < 1$, applying the BH95 procedure at level $\alpha/\pi_0$ would increase the number of rejections while keeping $\mathsf{FDR} \leq \alpha$. However, as $\pi_0$ is unknown, this is only an Oracle procedure. Many proposed procedures therefore try to imitate the Oracle by applying the BH95 procedure at level $\alpha/\widehat{\pi_0}$, where $\widehat{\pi_0}$ estimates (or at least provides an upper bound for) the true $\pi_0$ [2, 21, 23]. Such procedures are referred to as *two-stage adaptive* or *plug-in* procedures. A new class of procedures has recently been proposed, with the aim of providing tighter FDR control than the BH95 procedure, while avoiding explicit solution of the semi-parametric problem of $\pi_0$ estimation. Procedures of this second class are referred to as *one-stage adaptive* procedures [4, 7].

The FDR controlling properties of such procedures have been carefully studied for a finite number of hypotheses hypotheses [1, 2, 4, 8, 9, 19, 21], or asymptotically [5, 7, 8, 10, 11, 21, 22, 23]. As the proportion of erroneous rejections (FDP) is a stochastic quantity, its fluctuations around its mean value are worth investigating. Several procedures have been proposed for controlling the upper quantiles of the FDP [11, 12, 13, 15, 16, 17, 18, 26]. The asymptotic behavior of process $(\mathsf{FDP}_m(t))_{0<t\leq 1}$, where $t$ is a deterministic threshold, has also been studied [11, 21, 23]. We focus here on the properties of the *random threshold* $\widehat{\tau}$ associated with a given multiple testing procedure, particularly in the asymptotic distribution of $\mathsf{FDP}_m(\widehat{\tau})$, the FDP actually reached by the procedure.

**Organization of the paper.** In section 2 we propose a general framework for asymptotic analysis of the FDP of multiple testing procedures. In section 3 we derive the asymptotic distribution of the FDP of a multiple testing procedure with generic threshold function $\mathcal{T}$ and characterize the asymptotic equivalence of multiple testing procedures. These results are explicitly connected to the regularity of the map $\mathcal{T}$, which is then discussed. In section 4 we derive the asymptotic behavior of several existing procedures. In section 5 we point out interesting connections between one-stage adaptive and two-stage adaptive procedures. The main results are summarized and discussed in section 6, and proofs of the main results are gathered in section 7.

## 2. Background and notation

### *2.1. Background*

Throughout the paper, we consider a sequence $(P_i)_{i\in\mathbb{N}}$ of $p$-values associated with a collection of binary tests of a null hypothesis $\mathcal{H}_0$ against an alternative hypothesis $\mathcal{H}_1$.

**Definition 2.1** (Multiple Testing Procedure (MTP))**.** *A multiple testing procedure $\mathcal{M}$ is a sequence of functions $\mathcal{M}_m : [0,1]^m \to [0,1]$ such that for any $m$-dimensional vector of $p$-values $(P_1, \dots P_m)$, all hypotheses $i$ satisfying*

$$P_i \leq \mathcal{M}_m(P_1, \dots P_m)\,.$$

*are rejected. Slightly abusing notation, we shall write $\mathcal{M}(P_1, \dots P_m)$ for $\mathcal{M}_m(P_1, \dots P_m)$.*

Denoting by $V_m$ and $R_m$ the number of illegitimate rejections and the total number of rejections among the $m$ tested hypotheses for a multiple testing procedure $\mathcal{M}$, the associated False Discovery Proportion and False Discovery Rate are $\mathsf{FDP}_m(\mathcal{M}) = \frac{V_m}{R_m \vee 1}$, and $\mathsf{FDR}_m(\mathcal{M}) = \mathbb{E}\left[\frac{V_m}{R_m \vee 1}\right]$.

**Model.** We consider the setting originally proposed by [1]: among $m$ tested hypotheses, $m_0(m)$ are true nulls. We assume that $p$-values are uniformly distributed on $[0,1]$ under $\mathcal{H}_0$, and distributed according to $G_1$ under $\mathcal{H}_1$, where $G_1$ is a concave, $C^1$ distribution function, with density $g_1$.

We assume that $\pi_0(m) = m_0(m)/m$ verifies $\lim_{m\to+\infty} \pi_0(m) = \pi_0$, which we call *proportion of true null hypotheses*. We let $G(x) = \pi_0 x + (1-\pi_0)G_1(x)$, and $g = \pi_0 + (1-\pi_0)g_1$ be the corresponding density. Finally, we assume that *all $p$-values are independent.*

**The BH95 procedure.** Letting $(P_{(i)})_{1\leq i\leq m}$ be the vector of ordered $p$-values associated with $(P_i)_{1\leq i\leq m}$, and

$$\widehat{I}_m = \max\left\{i \in \{1, \dots m\}, P_{(i)} \leq \alpha i/m\right\}\,,$$

the BH95 procedure at level $\alpha$ is defined by $\mathcal{M}^\alpha_m(P_1, \dots P_m) = \alpha \widehat{I}_m/m$. This definition can be rewritten as follows. Letting $\widehat{\mathbb{G}}_m$ be the empirical distribution function of the $p$-values,

$$\mathcal{M}^\alpha(P_1, \dots P_m) = \sup\left\{u \in [0,1], \widehat{\mathbb{G}}_m(u) \geq u/\alpha\right\}\,.$$

These two equivalent formulations of the BH95 threshold are illustrated in Figure 1. Following [7], $u \mapsto u/\alpha$ will be called the *rejection curve* of the BH95 procedure (also known as *Simes' line* [20]).

## 2.2. Formalism

**Threshold functions.** This interpretation of the BH95 procedure in terms of the empirical distribution function suggests to define *threshold functions* as follows. Let $D[0,1]$ denote the set of cadlag functions defined on $[0,1]$, that is, the set of functions defined on $[0,1]$ that are right-continuous and have left limits.

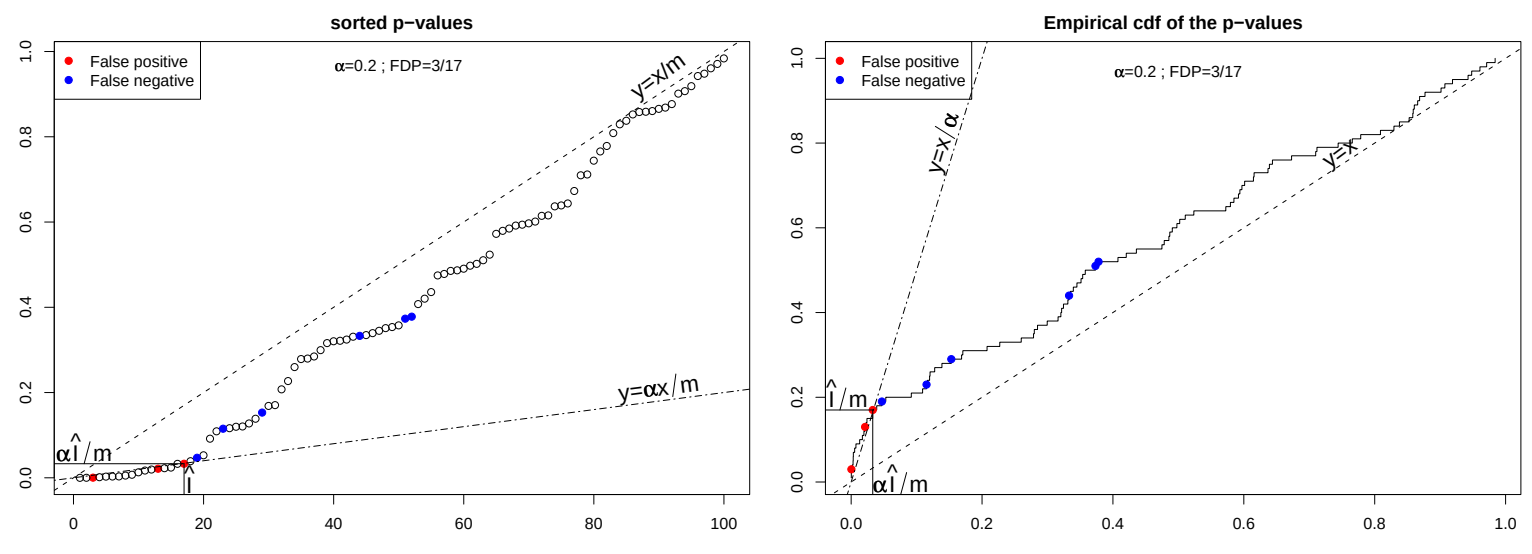

FIGURE 1. Dual interpretations of the BH95 threshold.

**Definition 2.2** (Threshold function)**.** *A multiple testing procedure* $\mathcal{M}$ *has threshold function* $\mathcal{T} : D[0,1] \to [0,1]$ *if and only if*

$$\forall m \in \mathbb{N}, \mathcal{M}(P_1, \dots P_m) = \mathcal{T}(\widehat{\mathbb{G}}_m)\,.$$

Note that $\mathcal{T}$ does not depend on $m$ in Definition 2.2. In the remainder of this paper $\mathcal{T}(G)$ will be denoted by $\tau^\star$.

**FDP as a stochastic process of a random threshold.** As suggested in a previous study [11], the False Discovery Proportion can be viewed as a stochastic process. Let $\widehat{\mathbb{G}}_{0,m}$ and $\widehat{\mathbb{G}}_{1,m}$ denote the (unobservable) empirical distribution function of the $p$-values under the null and alternative hypotheses:

$$\begin{cases} \widehat{\mathbb{G}}_{0,m}(t) & = \dfrac{\sum_{\{i/\mathcal{H}_0 \text{ true}\}} \mathbf{1}_{P_i \le t}}{m_0(m)} \\ \widehat{\mathbb{G}}_{1,m}(t) & = \dfrac{\sum_{\{i/\mathcal{H}_1 \text{ true}\}} \mathbf{1}_{P_i \le t}}{m - m_0(m)} \end{cases}.$$

As $\pi_0(m)$ is deterministic in our setting and verifies $\lim_{m\to+\infty} \pi_0(m) = \pi_0$, we will assume without loss of generality that it is constant equal to $\pi_0$, in order to alleviate notation. Therefore, we have $\widehat{\mathbb{G}}_m = \pi_0 \widehat{\mathbb{G}}_{0,m} + (1-\pi_0)\widehat{\mathbb{G}}_{1,m}$, and, for any $t \in [0,1]$, $R_m(t) = \frac{1}{m}\sum_{i=1}^m \mathbf{1}_{P_i \le t} = \widehat{\mathbb{G}}_m(t)$ and $V_m(t) = \frac{1}{m}\sum_{\{i/\mathcal{H}_0 \text{ true}\}} \mathbf{1}_{P_i \le t} = \pi_0 \widehat{\mathbb{G}}_{0,m}(t)$, so that

$$\mathsf{FDP}_m(t) = \frac{\pi_0 \widehat{\mathbb{G}}_{0,m}(t)}{\widehat{\mathbb{G}}_m(t) \vee \frac{1}{m}}$$

is the False Discovery Proportion achieved at the deterministic threshold $t$. The asymptotic properties of the stochastic process $(\mathsf{FDP}_m(t))_{0\le t\le 1}$ were analyzed by Genovese and Wasserman [11]. They noticed that $\mathsf{FDR}_m(t) = \mathbb{E}[\mathsf{FDP}_m(t)]$, so the achieved FDR at $t$, may be written as

$$\mathsf{FDR}_m(t) = \mathsf{p}(t)\left(1 - (1 - G(t))^m\right),$$

where $\mathsf{p}(t) = \frac{\pi_0 t}{G(t)}$ is the *positive False Discovery Rate* (pFDR) at $t$, as defined by [21]. They proved that the $\mathsf{FDP}_m$ process converges to pFDR at a rate $\frac{1}{\sqrt{m}}$, and built *confidence envelopes* for the FDP process using this result.

We make use of this stochastic process approach here to study the behavior of the FDP *actually achieved by a given multiple testing procedure* $\mathcal{T}$, that is, the random variable $\mathsf{FDP}_m(\mathcal{T}(\widehat{\mathbb{G}}_m))$. We investigated the asymptotic behavior of this variable and, in particular, its fluctuations around the asymptotic FDR achieved by procedure $\mathcal{T}$, by writing $\mathsf{FDP}_m(\mathcal{T}(\widehat{\mathbb{G}}_m))$ as a function of the empirical distribution functions under the null and alternative hypotheses. Letting

$$\mathcal{V} : (F_0, F_1) \mapsto \pi_0 F_0(\mathcal{T}(\pi_0 F_0 + (1 - \pi_0) F_1))$$

and

$$\mathcal{R} : F \mapsto F(\mathcal{T}(F)),$$

the FDP achieved by procedure $\mathcal{T}$ may be written as

$$\mathsf{FDP}_m(\mathcal{T}(\widehat{\mathbb{G}}_m)) = \frac{\mathcal{V}(\widehat{\mathbb{G}}_{0,m}, \widehat{\mathbb{G}}_{1,m})}{\mathcal{R}(\pi_0 \widehat{\mathbb{G}}_{0,m} + (1 - \pi_0)\widehat{\mathbb{G}}_{1,m}) \vee \frac{1}{m}}$$

since $\widehat{\mathbb{G}}_m = \pi_0 \widehat{\mathbb{G}}_{0,m} + (1-\pi_0)\widehat{\mathbb{G}}_{1,m}$. Using the functional Delta method [27], this formalism makes it possible to break down the analysis of $\mathsf{FDP}_m(\mathcal{T}(\widehat{\mathbb{G}}_m))$ into the regularity properties of the map $\mathcal{T}$, which depend solely on the procedure, and the asymptotic behavior of the empirical distribution functions of the $p$-values, which can be derived from Donsker's invariance principle [6] because $p$-values are assumed to be independent.

*Remark* 2.3. This paper focuses on FDP, but the formalism we propose here could be used to derive the asymptotic distribution of any risk measure based on the number of true/false positive/negatives, under the same regularity conditions. In particular, the results obtained here can also be applied to the False Non-discovery Proportion (FNP) [10]:

$$\mathsf{FNP}_m(t) = \frac{(1 - R_m(t)/m) - (\pi_0 - V_m(t)/m)}{1 - \pi_0}$$

**Multiple testing procedures studied.** The threshold function of the BH95 procedure is defined by

$$\mathcal{T}(F) = \sup\{u \in [0, 1], F(u) \geq u/\alpha\}\,.$$

As the BH95 procedure keeps the false discovery rate at a level of (exactly) $\pi_0\alpha$ when $p$-values are independent [3, 8, 19, 23], it is conservative by a factor $\pi_0$. Other multiple testing procedures have been proposed that estimate $\pi_0$, either implicitly or explicitly, to provide tighter (i.e. more powerful) FDR control under independence:

**One-stage adaptive procedures (**BR08 **[4],** FDR08 **[7])** use rejection curves other than Simes' line, without explicitly incorporating an estimate of $\pi_0$.

**Two-stage adaptive procedures (**BKY06 **[2],** STS04 **[23],** Sto02 **[21])** apply the BH95 procedure at a level of $\alpha/\widehat{\pi_0}$, where $\widehat{\pi_0}$ is an estimator of $\pi_0$.

We therefore consider threshold functions of the form

$$\mathcal{T}(F) = \mathcal{U}\left(F, \mathcal{A}(F)\right),$$

with

$$\mathcal{U}(F,\alpha) = \sup\{u \in [0,1], F(u) \geq r_\alpha(u)\},$$

where $r_\alpha : [0,1] \to \mathbb{R}_+$ will be called a *rejection curve* (after [7]), and $\mathcal{A} : D[0,1] \to [0,1]$ will be called a *level function.* $r_\alpha$ will be denoted by $r(\alpha,\cdot)$ whenever the dependence on $\alpha$ is of importance. $\mathcal{A}$ and $r_\alpha$ are two degrees of freedom that can be used to describe generalizations of the BH95 procedure, corresponding to the case in which the level function is constant (equal to $\alpha$), *and* the rejection curve is Simes' line. In this paper we consider increasing rejection curves satisfying $r_\alpha(0) = 0$, so that $\mathcal{U}(F,\alpha) \geq 0$ for any $F \in D[0,1]$ and $\alpha \in [0,1]$.

### 2.3. Overview of main results

Theorem 3.2 below shows that the FDP of a multiple testing procedure with threshold function $\mathcal{T}$ converges in distribution at rate $1/\sqrt{m}$ to a conservative, procedure-specific FDR level. This theorem holds under a general regularity condition on the map $\mathcal{T}$, which is implied by the existence and uniqueness of an interior right-crossing point between the distribution function of the $p$-values and the rejection curve of the procedure; the existence condition for a given procedure may be interpreted as a natural generalization of the notion of *criticality*, which has recently been introduced for the BH95 procedure [5].

Although the BH95 procedure is known to control FDR at a level of *exactly* $\pi_0\alpha$ [3, 8], other procedures have been proved to yield only an FDR *not larger than* $\alpha$, either for a finite number of hypotheses (procedures STS04, BKY06 and BR08) or asymptotically (Sto02 and FDR08). In section 4 we derive the asymptotic behavior of these procedures, and the associated regularity conditions. As all procedures converge at the same rate $1/\sqrt{m}$, their asymptotic power may be explicitly compared through their attained asymptotic FDR.

In section 5 we demonstrate the existence of interesting connections between the one-stage and two-stage adaptive procedures under investigation: with a striking symmetry, procedure BR08 may be interpreted as a fixed point of the iteration of procedure BKY06, and procedure FDR08 as a fixed point of the iteration of procedure Sto02.

## 3. Asymptotic properties of threshold procedures

This section provides general results about multiple testing procedures with threshold functions satisfying the following regularity condition. We refer to [27] for a formal definition of Hadamard differentiability.

**Condition C.1** (Hadamard-differentiability)**.** *The threshold function $\mathcal{T}$ satisfies $\mathcal{T}(G) > 0$, and is Hadamard-differentiable at $G$, tangentially to $C[0,1]$, where $C[0,1]$ is the set of continuous functions on $[0,1]$ The threshold function derivative is denoted by $\dot{\mathcal{T}}_G$.*

We begin by deriving the asymptotic distribution of the FDP of any multiple testing procedure satisfying Condition C.1 (section 3.1). We then define and characterize asymptotic equivalence between multiple testing procedures in terms of Condition C.1 (section 3.2). Finally we interpret this Condition in terms of crossing points between the distribution function $G$ of the $p$-values and the rejection curve (section 3.3).

### 3.1. Asymptotic False Discovery Proportion

Condition C.1 makes it possible to use the functional Delta method [27] to derive the asymptotic distribution of the False Discovery Proportion $\mathsf{FDP}_m(\mathcal{T}(\widehat{\mathbb{G}}_m))$ actually achieved by procedure $\mathcal{T}$ from the convergence in distribution of the centered empirical processes associated with $\widehat{\mathbb{G}}_{0,m}$ and $\widehat{\mathbb{G}}_{1,m}$, which is a consequence of Donsker's theorem [27]:

**Theorem 3.1** (Donsker)**.** *If the $p$-values are independent, then*

*(i)* $\sqrt{m}\left(\begin{pmatrix}\widehat{\mathbb{G}}_{0,m}\\ \widehat{\mathbb{G}}_{1,m}\end{pmatrix} - \begin{pmatrix}G_0\\ G_1\end{pmatrix}\right) \rightsquigarrow \begin{pmatrix}\mathbb{Z}_0\\ \mathbb{Z}_1\end{pmatrix}$ *on $[0,1]$, where $\mathbb{Z}_0$ and $\mathbb{Z}_1$ are independent Gaussian processes such that $\mathbb{Z}_0 \overset{(d)}{=} \mathbb{B}$ and $\mathbb{Z}_1 \overset{(d)}{=} \mathbb{B} \circ G_1$, and $\mathbb{B}$ is a standard Brownian bridge on $[0,1]$.*

*(ii)* $\sqrt{m}\big(\widehat{\mathbb{G}}_m - G\big) \rightsquigarrow \mathbb{Z}$ *on $[0,1]$, where $\mathbb{Z} = \pi_0\mathbb{Z}_0 + (1-\pi_0)\mathbb{Z}_1$ is a Gaussian process with continuous sample paths and covariance function given by*

$$\mathbb{E}\left[\mathbb{Z}(s)\mathbb{Z}(t)\right] \quad = \quad \pi_0^2\gamma_0(s,t) + (1-\pi_0)^2\,\gamma_0(G_1(s), G_1(t)),$$

*where $\gamma_0$ is the covariance function of $\mathbb{B}$, that is, $\gamma_0 : (s,t) \mapsto s \wedge t - st$.*

**Theorem 3.2** (Asymptotic distribution of $\mathsf{FDP}_m$ for procedure $\mathcal{T}$)**.** *Let $\mathcal{T}$ be a threshold function, $\tau^\star = \mathcal{T}(G)$, and $\mathsf{p}(t) = \frac{\pi_0 t}{G(t)}$ the positive False Discovery Rate at threshold $t$. Under Condition C.1,*

*(i)*

$$\sqrt{m}\left(\mathcal{T}(\widehat{\mathbb{G}}_m) - \tau^\star)\right) \rightsquigarrow \dot{\mathcal{T}}_G(\mathbb{Z})\,,$$

*(ii)*

$$\lim_{m\to\infty} \mathsf{FDR}_m(\mathcal{T}(\widehat{\mathbb{G}}_m)) = \mathsf{p}(\tau^\star)\,,$$

*(iii)*

$$\sqrt{m}\left(\mathsf{FDP}_m(\mathcal{T}(\widehat{\mathbb{G}}_m)) - \mathsf{p}(\tau^\star)\right) \rightsquigarrow X\,,$$

*with*

$$X = \mathsf{p}(\tau^\star)(1-\mathsf{p}(\tau^\star))\left(\frac{\mathbb{Z}_0(\tau^\star)}{\tau^\star} - \frac{\mathbb{Z}_1(\tau^\star)}{G_1(\tau^\star)}\right) + \dot{\mathsf{p}}(\tau^\star)\dot{\mathcal{T}}_G(\mathbb{Z})$$

*and* $\mathbb{Z} = \pi_0\mathbb{Z}_0 + (1-\pi_0)\mathbb{Z}_1$*, where* $\mathbb{Z}_0$ *and* $\mathbb{Z}_1$ *are independent Gaussian processes such that* $\mathbb{Z}_0 \overset{(d)}{=} \mathbb{B}$ *and* $\mathbb{Z}_1 \overset{(d)}{=} \mathbb{B}\circ G_1$*, and* $\mathbb{B}$ *is a standard Brownian bridge on* $[0,1]$*.*

According to $(ii)$, the asymptotic FDR achieved by procedure $\mathcal{T}$ is the pFDR at the asymptotic threshold $\tau^\star = \mathcal{T}(G)$. This is true because $\tau^\star$ is positive (by Condition C.1). In particular, Theorem 3.2 provides a necessary and sufficient condition under which a multiple testing procedure with Hadamard differentiable threshold function asymptotically controls FDR:

**Corollary 3.3.** *A threshold function* $\mathcal{T}$ *satisfying Condition C.1 asymptotically controls* FDR *if and only if its* pFDR *at* $\tau^\star = \mathcal{T}(G)$ *(i.e. its asymptotic* FDR*) is below* $\alpha$*, that is, if and only if*

$$\frac{\pi_0\tau^\star}{G(\tau^\star)} \leq \alpha\,.$$

*Remark* 3.4 (Form of $\dot{\mathcal{T}}_G$)*.* The expression of $\dot{\mathcal{T}}_G$ for threshold functions is given by Corollary 7.12, which shows that for one-stage adaptive procedures (where the level function $\mathcal{A}$ is constant), $\dot{\mathcal{T}}_G$ is proportional to the inverse of the difference between the slopes of $r_\alpha$ and $G$ at $\tau^\star$. For two-stage plug-in procedures, which typically estimate $\pi_0$ using $G(u_0)$ for some $u_0$ (e.g. $u_0 = \lambda$ for procedure Sto02), $\dot{\mathcal{T}}_G$ involves an additional term that depends on $G(u_0)$, and the asymptotic distribution of the FDP depends on the centered Gaussian random variable $\mathbb{Z}(u_0)$, where $\mathbb{Z}$ is defined in Theorem 3.1.

### *3.2. Asymptotically equivalent procedures*

Some multiple testing procedures cannot be written in terms of threshold functions, because they do not depend exclusively on $\widehat{\mathbb{G}}_m$, but instead also directly depend on the number $m$ of observations. When such procedures are only slight perturbations of actual threshold procedures, they share the same asymptotic distribution, as explained below.

**Definition 3.5** (Asymptotic equivalence of multiple testing procedures)**.** *Let* $\mathcal{T}$ *be a threshold function for which Condition C.1 holds for* $\mathcal{T}$*. A multiple testing procedure* $\mathcal{M}$ *is asymptotically equivalent to* $\mathcal{T}$ *as* $m\to+\infty$ *if and only if*

$$\sqrt{m}\left(\mathsf{FDP}_m\left(\mathcal{M}(P_1,\dots P_m)\right) - \mathsf{FDP}_m\left(\mathcal{T}(\widehat{\mathbb{G}}_m)\right)\right) \overset{P}{\to} 0\,.$$

**Proposition 3.6** (Asymptotic equivalence of thresholding procedures)**.** *Let* $\mathcal{T}$ *be a threshold function, and* $\varepsilon = (\varepsilon_m)_{m\in\mathbb{N}}$ *a positive sequence. For* $m\in\mathbb{N}$*, let* $\mathcal{T}_m : D[0,1]\to[0,1]$ *such that*

$$\forall F\in D[0,1], \mathcal{T}(F-\varepsilon_m) \leq \mathcal{T}_m(F) \leq \mathcal{T}(F)\,.$$

*If Condition C.1 holds for $\mathcal{T}$, and if $\varepsilon_m = o\left(\frac{1}{\sqrt{m}}\right)$, $\mathcal{T}_m$ is asymptotically equivalent to $\mathcal{T}$ as $m \to +\infty$.*

Several applications of Proposition 3.6 are given in section 4. For example, the asymptotic behavior of procedure $\mathcal{T}_m = \mathsf{STS04}(\lambda)$ can be derived from that of procedure $\mathcal{T} = \mathsf{Sto02}(\lambda)$, for which Theorem 3.2 may be used because $\mathsf{Sto02}(\lambda)$ is an actual threshold function.

### *3.3. Regularity conditions*

For the threshold functions under investigation, $\mathcal{T}(G)$ is defined as the last point for which $G \geq r(\mathcal{A}(G), \cdot)$. Therefore, the existence of a *unique interior right crossing point* between $G$ and $r(\mathcal{A}(G), \cdot)$ ensures that Theorem 3.2 and Proposition 3.6 are applicable, i.e. that $\mathcal{T}(G) > 0$, and that $\mathcal{T}$ is Hadamard differentiable at $G$ (Condition C.1). For two-stage adaptive (plug-in) procedures, for which the level function $\mathcal{A}$ is not constant, additional technical assumptions concerning the regularity of $\mathcal{A}$ require checking (see Corollary 7.12) to ensure that Condition C.1 holds.

**Definition 3.7** (Right crossing point)**.** *Let $r_\alpha$ be a rejection curve, and $\mathcal{A}$ a level function. Denote by $\mathcal{T} : F \mapsto \mathcal{U}(F, \mathcal{A}(F))$ the associated threshold function, where $\mathcal{U}(F, \alpha) = \sup\{u \in [0,1], F(u) \geq r_\alpha(u)\}$. A right crossing point for the multiple comparison problem defined by $\mathcal{T}$ (or, in short, a right crossing point for $\mathcal{T}$), is a point $t \in [0,1]$ such that $G(t) = r_\alpha(t)$, and $g(t) < \frac{\partial r}{\partial u}(\mathcal{A}(G), t)$. If $t$ belongs to the open interval $(0,1)$ it is called an interior right crossing point for $\mathcal{T}$.*

Condition $g(t) < \frac{\partial r}{\partial u}(\mathcal{A}(G), t)$ in Definition 3.7 ensures that $G$ and $r_{\mathcal{A}(G)} = r(\mathcal{A}(G), \cdot)$ actually *cross* at $t$, *i.e.* that $G \geq r_{\mathcal{A}(G)}$ in a left-neighborhood of $t$, and that $G \leq r_{\mathcal{A}(G)}$ in a right-neighborhood of $t$.

Studies of the asymptotic distribution of the abovementioned FDR controlling procedures require investigation, in each case, of the conditions guaranteeing the existence of a unique interior right crossing point. To this end, we broke this condition down as follows:

**Condition C.2** (Existence)**.** *$\mathcal{T}$ has an interior right crossing point.*

**Condition C.3** (Uniqueness)**.** *$\mathcal{T}$ has at most one interior right crossing point.*

Condition C.3 always holds for procedures based on Simes' line (BH95, Sto02, and BKY06) because their rejection curve is linear, and $G$ is concave. Condition C.2 typically holds in situations in which the slope of $G$ at the origin is large enough. In the case of the BH95 procedure, Chi recently showed the existence of a *critical value* $\alpha^\star$ depending solely on the distribution function $G$ of the $p$-values, such that if $\alpha < \alpha^\star$, the number of discoveries made by the BH95 procedure is stochastically bounded as the number of tested hypotheses increases, whereas if $\alpha > \alpha^\star$, the proportion of discoveries converges in probability to a positive value $\tau^\star = \mathcal{T}(G)$ [5].

In section 4, we provide a detailed analysis of a number of FDR controlling procedures, and present, for each, a *critical value* for the target FDR level characterising situations in which condition C.2 is guaranteed for the procedure.

## 4. Results for procedures of interest

We apply the results of the preceding section to a series of procedures with proven (asymptotic) FDR control. Starting from the original BH95 procedure and its Oracle version (section 4.1), we then turn to adaptive procedures, which implicitly or explicitly incorporate an estimate of the proportion $\pi_0$ of true null hypotheses: *one-stage adaptive procedures* are studied in section 4.2, and *two-stage adaptive procedures* (also called plug-in procedures) are studied in section 4.3.

### *4.1.* BH95 *procedure*

We will first recall the definition of the BH95 procedure in our framework.

**Definition 4.1** (Procedure BH95[1])**.** *The* BH95 *procedure is the multiple testing procedure with threshold function*

$$\mathcal{T}^{\mathsf{BH95}}(F) = \sup\{u \in [0,1], F(u) \geq u/\alpha\}\,.$$

As the rejection curve of procedure BH95 is linear, and $G$ is concave, the uniqueness Condition C.3 always holds, and the existence Condition C.2 can be reduced to $\alpha > \alpha^\star$, where $\alpha^\star = \inf_{u\to 0} u/G(u) = \lim_{u\to 0} 1/g(u)$ corresponds to the *critical value of the* BH95 *procedure* [5]:

**Condition C.4** (Condition C.2 for the BH95 procedure)**.** *The target* FDR *level* $\alpha$ *is greater than the critical value* $\alpha^\star$ *of the* BH95 *procedure.*

The criticality phenomenon is illustrated in Figure 2 for Laplace (double exponential) test statistics. The Weak Law of Large Numbers phenomenon analyzed by [5], which occurs when $\alpha > \alpha^\star$, was noted by [10]. We now derive the corresponding central limit theorem under the same hypothesis, and the asymptotic distribution of the FDP actually achieved by the BH95 procedure.

**Theorem 4.2** (Asymptotic properties of the BH95 procedure)**.** *Let* $\tau^\star = \mathcal{T}^{\mathsf{BH95}}(G)$*. Under Condition C.4,*

*(i)*

$$\sqrt{m}\left(\mathcal{T}^{\mathsf{BH95}}(\widehat{\mathbb{G}}_m) - \tau^\star\right) \rightsquigarrow \frac{\mathbb{Z}(\tau^\star)}{1/\alpha - g(\tau^\star)}\,,$$

*with* $\mathbb{Z} = \pi_0\mathbb{Z}_0 + (1-\pi_0)\mathbb{Z}_1$*, where* $\mathbb{Z}_0$ *and* $\mathbb{Z}_1$ *are independent Gaussian processes such that* $\mathbb{Z}_0 \overset{(d)}{=} \mathbb{B}$ *and* $\mathbb{Z}_1 \overset{(d)}{=} \mathbb{B}\circ G_1$*, and* $\mathbb{B}$ *is a standard Brownian bridge on* $[0,1]$*.*

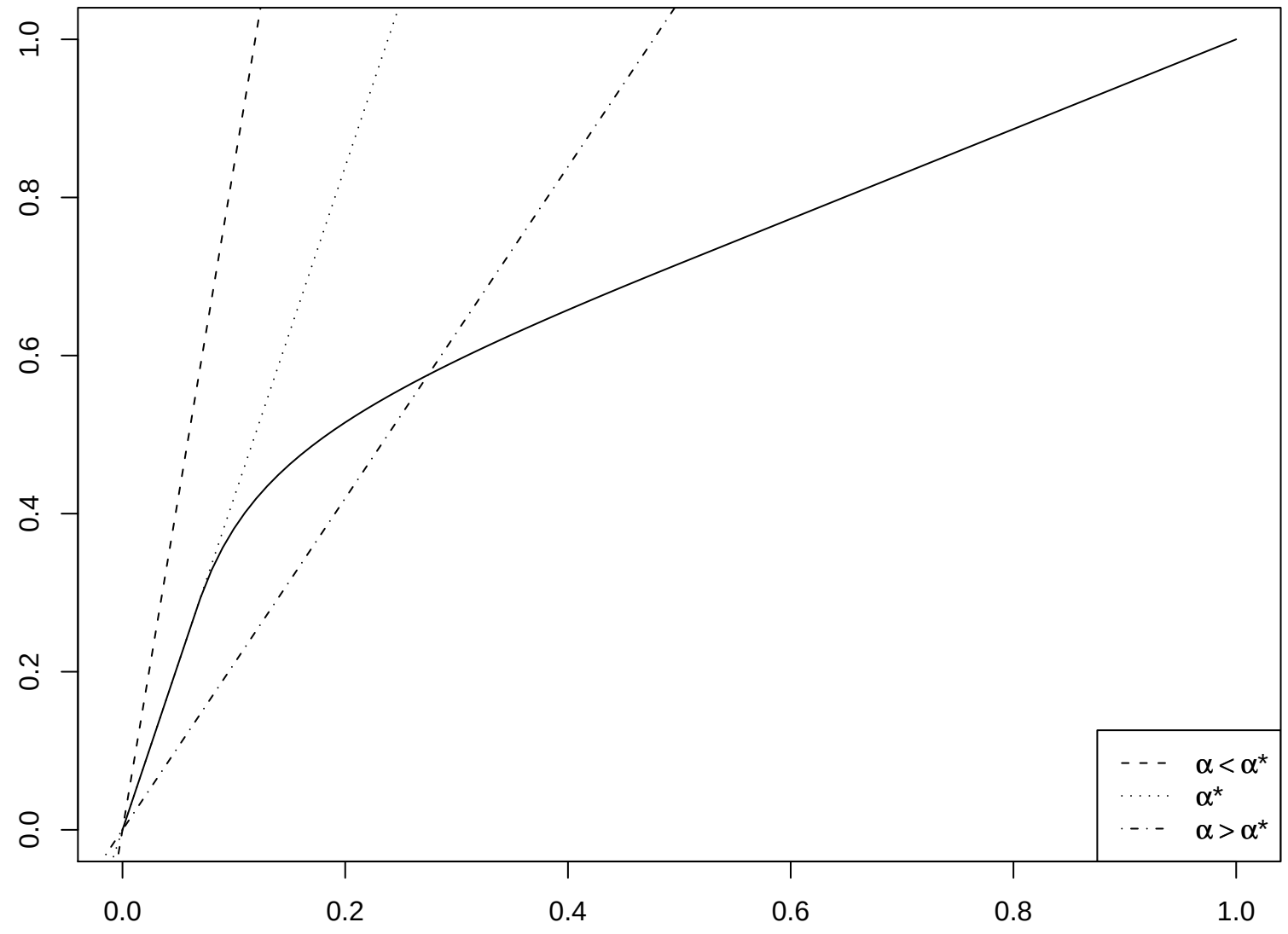

FIGURE 2. Critical value of the BH95 procedure for Laplace test statistics with location parameter $\theta = 2$, and $\pi_0 = 0.5$. Solid line: distribution function $G$; straight lines: Simes' lines for several values of $\alpha$. There is an interior right crossing point between the distribution function of the $p$ values and Simes' line if and only if $\alpha > \alpha^\star = 1/g(0)$.

*(ii)*

$$\sqrt{m}\left(\mathsf{FDP}_m(\mathcal{T}^{\mathsf{BH95}}(\widehat{\mathbb{G}}_m)) - \pi_0\alpha\right) \rightsquigarrow \mathcal{N}\left(0, (\pi_0\alpha)^2\frac{1-\tau^\star}{\tau^\star}\right).$$

Applying the BH95 procedure at level $\alpha/\pi_0$ leads to an Oracle procedure (as $\pi_0$ is not known) that is more powerful as it controls FDR at level exactly $\alpha$. This procedure, which we denote by BH95o, has threshold function

$$\mathcal{T}^{\mathsf{BH95o}}(F) = \sup\{u \in [0,1], F(u) \geq \pi_0 u/\alpha\},$$

and its critical value is therefore $\pi_0\alpha^\star$, which translates into the following regularity condition:

**Condition C.5** (Condition C.2 for the BH95 Oracle procedure)**.** *The target* FDR *level $\alpha$ is greater than $\pi_0\alpha^\star$, where $\alpha^\star$ is the critical value of the* BH95 *procedure.*

The corresponding asymptotic properties can be derived from Theorem 4.2:

**Corollary 4.3** (Asymptotic properties of the BH95 Oracle procedure)**.** *Let $\tau^\star = \mathcal{T}^{\mathsf{BH95o}}(G)$. Under Condition C.5,*

*(i)*

$$\sqrt{m}\left(\mathcal{T}^{\mathsf{BH95o}}(\widehat{\mathbb{G}}_m) - \tau^\star\right) \rightsquigarrow \frac{\mathbb{Z}(\tau^\star)}{\pi_0/\alpha - g(\tau^\star)},$$

*with* $\mathbb{Z} = \pi_0\mathbb{Z}_0 + (1-\pi_0)\mathbb{Z}_1$*, where* $\mathbb{Z}_0$ *and* $\mathbb{Z}_1$ *are independent Gaussian processes such that* $\mathbb{Z}_0 \overset{(d)}{=} \mathbb{B}$ *and* $\mathbb{Z}_1 \overset{(d)}{=} \mathbb{B}\circ G_1$*, and* $\mathbb{B}$ *is a standard Brownian bridge on* $[0,1]$*.*

*(ii)*

$$\sqrt{m}\left(\mathsf{FDP}_m(\mathcal{T}^{\mathsf{BH95o}}(\widehat{\mathbb{G}}_m)) - \alpha\right) \rightsquigarrow \mathcal{N}\left(0, \alpha^2\frac{1-\tau^\star}{\tau^\star}\right).$$

### *4.2. One-stage adaptive procedures*

The first class of adaptive procedures studied here are *one-stage* adaptive procedures, because they estimate $\pi_0$ implicitly, rather than through a level function $\mathcal{A}$.

**Definition 4.4** (Adaptive procedure)**.** *Let* $r_\alpha : [0,1] \to [0,1]$*. The adaptive procedure associated with* $r_\alpha$ *is the multiple testing procedure defined by the threshold function*

$$\mathcal{T}(F) = \sup\{u \in [0,1], F(u) \geq r_\alpha(u)\}.$$

The rejection curve of adaptive procedures is not linear, so the conditions under which Condition C.1 is fulfilled are more subtle than for the BH95 procedure (section 4.1) or for two-stage adaptive procedures (section 4.3).

**Procedure** FDR08($\lambda$)**.** The rejection curve of the FDR08 procedure [7] is defined for $u \in [0,1]$ by $f_\alpha(u) = \frac{u}{\alpha+(1-\alpha)u}$. As $f_\alpha(1) = 1$, the corresponding threshold function is always equal to 1. This procedure therefore systematically rejects *all* hypotheses, and does not control FDR either for finite sample size or asymptotically. Several ways of overcoming this problem have been proposed [7], including *truncating* the rejection curve, yielding the following procedure:

**Definition 4.5** (Procedure FDR08($\lambda$))**.** *Let* $\lambda \in [0,1)$*. The rejection curve of the* FDR08($\lambda$) *procedure is defined by* $f_\alpha^\lambda(u) = f_\alpha(u)$ *for* $u \leq \lambda$*, and* $+\infty$ *otherwise. The threshold function of the* FDR08($\lambda$) *procedure is therefore given by*

$$\mathcal{T}^{\mathsf{FDR08}}(F) = \sup\left\{u \in [0,\lambda], F(u) \geq \frac{u}{\alpha+(1-\alpha)u}\right\}.$$

We introduce the following regularity condition:

**Condition C.6.** $\lambda \geq \kappa$*, where* $\kappa = \frac{\alpha(1-\pi_0)}{(1-\alpha)\pi_0}$*.*

Note that $\left(\kappa, \frac{1-\pi_0}{1-\alpha}\right)$ is the crossing point between the rejection curve $f_\alpha$ and the distribution function $\mathrm{DU}(\pi_0)$ in the extremal Dirac-Uniform configuration where all $p$-values drawn from $\mathcal{H}_1$ are equal to 0. As $G = \pi_0 G_0 + (1-\pi_0)G_1 \leq \mathrm{DU}(\pi_0)$, condition C.6 ensures that any interior right crossing point between $G$ and $f_\alpha$ occurs before $\lambda$. In practice, $\kappa$ is unknown because it depends on $\pi_0$. However, an upper bound for $\kappa$ can be deduced from a lower bound for $\pi_0$; for example, in microarray data analysis, it can often be assumed that $\pi_0 > \frac{1}{2}$: in this case, $\kappa$ is smaller than $\frac{\alpha}{1-\alpha}$.

By definition 4.5, the rejection curve $f_\alpha^\lambda$ of any procedure $\mathsf{FDR08}(\lambda)$ satisfying Condition C.6 is equal to $f_\alpha$ on $[0, \kappa]$, corresponding to the admissible region for interior right crossing points. The following Proposition is a straightforward consequence of this observation:

**Proposition 4.6.** *All* $\mathsf{FDR08}(\lambda)$ *procedures satisfying Condition C.6 are asymptotically equivalent in the sense of Definition 3.5.*

As the corresponding asymptotic distribution does not depend on $\lambda$, we will refer to it simply as the "asymptotic distribution of the FDR08 procedure". In order to characterize this distribution we introduce a further technical condition to ensure that $\kappa < 1$. Combined with Condition C.4, it also ensures that existence Condition C.2 holds for procedure $\mathsf{FDR08}(\lambda)$, because the slope of $f_\alpha^\lambda$ at the origin is $1/\alpha$.

**Condition C.7.** $\alpha < \pi_0$.

Condition C.7 is a mild assumption in practice, because $\pi_0$ is typically expected to be greater than $1/2$, in microarray data analysis, for example. When $\alpha \geq \pi_0$, there is no need for sophisticated FDR controlling procedures because rejecting all hypotheses yields $\mathsf{FDP} = \pi_0$ and thus $\mathsf{FDR} \leq \alpha$.

**Theorem 4.7** (Asymptotic behavior of procedure FDR08)**.** *Let* $\lambda \in [0,1)$ *such that Condition C.6 is fulfilled, and* $\tau^\star = \sup\left\{u \in [0,\kappa], G(u) \geq \frac{u}{\alpha+(1-\alpha)u}\right\}$. *Under uniqueness Condition C.3, and existence Conditions C.4 and C.7, we have*

$$\sqrt{m}\left(\mathsf{FDP}_m(\mathcal{T}^{\mathsf{FDR08}(\lambda)}(\widehat{\mathbb{G}}_m)) - \alpha\frac{\pi_0}{\overline{\pi_0}\,(\tau^\star)}\right) \rightsquigarrow X^{\mathsf{FDR08}},$$

*with* $\overline{\pi_0}\,(\tau^\star) = \frac{1-G(\tau^\star)}{1-\tau^\star}$, *and*

$$X^{\mathsf{FDR08}} = \mathsf{p}^\star(1-\mathsf{p}^\star\zeta(\tau^\star))\frac{\mathbb{Z}_0(\tau^\star)}{\tau^\star} - \mathsf{p}^\star(1-\mathsf{p}^\star)\zeta(\tau^\star)\frac{\mathbb{Z}_1(\tau^\star)}{G_1(\tau^\star)},$$

*where* $\mathsf{p}^\star = \alpha\pi_0/\overline{\pi_0}\,(\tau^\star)$ *is the* pFDR *achieved by procedure* FDR08*,*

$$\zeta(\tau^\star) = -\frac{(1-\overline{\pi_0}(\tau^\star))\overline{\pi_0}(\tau^\star)/\alpha}{\overline{\pi_0}(\tau^\star)^2/\alpha - g(\tau^\star)},$$

$\mathbb{Z}_0$ *and* $\mathbb{Z}_1$ *are independent Gaussian processes such that* $\mathbb{Z}_0 \stackrel{(d)}{=} \mathbb{B}$ *and* $\mathbb{Z}_1 \stackrel{(d)}{=} \mathbb{B}\circ G_1$, *and* $\mathbb{B}$ *is a standard Brownian bridge on* $[0,1]$.

As $\overline{\pi_0}\,(\tau^\star) = \frac{1-\tau^\star}{1-G(\tau^\star)} \in [\pi_0, 1]$, we have $\pi_0\alpha \leq \mathsf{p}^\star \leq \alpha$, so that procedure FDR08 is asymptotically more powerful than procedure BH95, and less powerful than procedure BH95o.

**Procedure** $\mathsf{BR08}(\lambda)$**.**

**Definition 4.8** (Procedure BR08($\lambda$) [4]). *Let $\lambda \in [0,1)$. The rejection curve of the* BR08($\lambda$) *procedure is defined by $b_\alpha^\lambda(u) = \frac{u}{\alpha(1-\lambda)+u}$ for $u \leq \lambda$, and $+\infty$ otherwise. The threshold function of the* BR08($\lambda$) *procedure is therefore given by*

$$\mathcal{T}^{\mathsf{BR08}(\lambda)}(F) = \sup\left\{u \in [0,\lambda], F(u) \geq \frac{u}{\alpha(1-\lambda)+u}\right\}.$$

Procedure BR08($\lambda$) is actually defined by the rejection curve $\left(1+\frac{1}{m}\right) b_\alpha^\lambda$ [4]. However these procedures are asymptotically equivalent according to Proposition 3.6; we will therefore use Definition 4.8.

As for the FDR08 procedure, the rejection curve of the BR08($\lambda$) procedure is not linear and we therefore need to make two assumptions to ensure that existence Condition C.2 holds: Condition C.8 ensures that there is no criticality phenomenon, that is, that the slope of the distribution function $G$ is great enough at the origin, and Condition C.9 ensures that a right crossing point occurs before $\lambda$, because the BR08($\lambda$) procedure is truncated at $\lambda$:

**Condition C.8.** *The target* FDR *level $\alpha$ satisfies $\alpha(1-\lambda) > \alpha^\star$, where $\alpha^\star$ is the critical value of the* BH95 *procedure.*

**Condition C.9.** *The distribution function $G$ satisfies*

$$G(\lambda) \leq \frac{\lambda}{\alpha}\frac{1-G(\lambda)}{1-\lambda}.$$

*Remark* 4.9. Condition C.9 may be written as $G(\lambda) \leq b_\alpha^\lambda$, or as $G(\lambda) \leq f_\alpha^\lambda$, because the rejection curves of procedures BR08($\lambda$) and FDR08 intersect at $\lambda$.

**Theorem 4.10** (Asymptotic distribution of procedure BR08($\lambda$)). *Let $\lambda \in [0,1)$ and $\tau^\star = \mathcal{T}^{\mathsf{BR08}(\lambda)}(G)$. Under uniqueness Conditions C.3 and existence Conditions C.8 and C.9, we have*

$$\sqrt{m}\left(\mathsf{FDP}_m(\mathcal{T}^{\mathsf{BR08}(\lambda)}(\widehat{\mathbb{G}}_m)) - \alpha\pi_0\frac{1-\lambda}{1-G(\tau^\star)}\right) \rightsquigarrow X^{\mathsf{BR08}(\lambda)},$$

*with*

$$X^{\mathsf{BR08}(\lambda)} = \mathsf{p}^\star(1-\mathsf{p}^\star\zeta(\tau^\star))\frac{\mathbb{Z}_0(\tau^\star)}{\tau^\star} - \mathsf{p}^\star(1-\mathsf{p}^\star)\zeta(\tau^\star)\frac{\mathbb{Z}_1(\tau^\star)}{G_1(\tau^\star)},$$

*where $\mathsf{p}^\star = \alpha\pi_0\frac{1-\lambda}{1-G(\tau^\star)}$ is the* pFDR *achieved by procedure* BR08($\lambda$)*,*

$$\zeta(\tau^\star) = -\frac{G(\tau^\star)^2/\tau^\star}{G(\tau^\star)(1-G(\tau^\star))/\tau^\star - g(\tau^\star)},$$

*$\mathbb{Z}_0$ and $\mathbb{Z}_1$ are independent Gaussian processes such that $\mathbb{Z}_0 \overset{(d)}{=} \mathbb{B}$ and $\mathbb{Z}_1 \overset{(d)}{=} \mathbb{B}\circ G_1$, and $\mathbb{B}$ is a standard Brownian bridge on $[0,1]$.*

Theorem 4.10 implies that procedure BR08($\lambda$) controls FDR asymptotically at level $\alpha$: as $\tau^\star \leq \lambda$, we have $\mathsf{p}^\star \leq \alpha\pi_0\frac{1-\lambda}{1-G(\lambda)}$, which is smaller than $\alpha$ because $\overline{\pi_0}(\lambda) = \frac{1-G(\lambda)}{1-\lambda}$ is an upper bound for $\pi_0$.

However, as $b_\alpha^\lambda(u) \geq u/\alpha$ if and only if $u \geq \lambda\alpha$, procedure $\mathsf{BR08}(\lambda)$ need not be more powerful than procedure BH95, and we have the following characterization:

$$\mathsf{BR08}(\lambda) \gg \mathsf{BH95} \iff \tau^\star_{\mathsf{BH95}} \geq \alpha\lambda\,,$$

where $\gg$ means "is more powerful than", and $\tau^\star_{\mathsf{BH95}}$ is the asymptotic threshold of procedure BH95. An explicit characterization of situations in which $\mathsf{BR08}(\lambda) \gg$ BH95 for Gaussian test statistics is given in [4].

### *4.3. Two-stage adaptive (plug-in) procedures*

In this section we study two-stage adaptive or *plug-in* procedures, in which a conservative step-up procedure is applied to a data-dependent level. In particular, we consider the case of *Simes' line-based* plug-in procedures, in which procedure BH95 is applied at level $\alpha/\widehat{\pi_0}$, where $\widehat{\pi_0}$ is estimated from the data:

**Definition 4.11** (Simes' line-based plug-in procedure)**.** *Let $\mathcal{A} : D[0,1] \to \mathbb{R}^*_+$. The Simes' line-based plug-in procedure associated with $\mathcal{A}$ is the multiple testing procedure defined by the threshold function*

$$\mathcal{T}(F) = \sup\left\{u \in [0,1], F(u) \geq \frac{u}{\mathcal{A}(F)}\right\}.$$

*Such procedures will simply be called* plug-in procedures *hereafter.*

As $r_\alpha$ is linear, and $G$ is concave, uniqueness Condition C.3 always holds for plug-in procedures, and existence Condition C.2 is the same as for procedure BH95, except that $\alpha$ is replaced by the value of the level function $\mathcal{A}$ at $G$:

**Condition C.10** (Condition C.2 for plug-in procedures)**.** *The level function $\mathcal{A}(G)$ associated with the target* FDR *level $\alpha$ is greater than the critical value of the* BH95 *procedure.*

Care is required when deriving the asymptotic distribution of the FDP for plug-in procedures, because the Hadamard derivative of $\mathcal{T}$, $\dot{\mathcal{T}}_G(H)$, typically involves the value of $H$ at $\tau^\star$ and at a point $u(\lambda)$ used for the estimation of $\pi_0$: $u(\lambda) = \lambda$ for procedure Sto02, and $u(\lambda) = \mathcal{U}(G,\lambda)$ for procedure $\mathsf{BKY06}(\lambda)$. The asymptotic variance of the False Discovery Proportion therefore involves the covariance between $\mathbb{Z}(\tau^\star)$ and $\mathbb{Z}(u(\lambda))$.

**Theorem 4.12** (Asymptotic FDP for procedures based on Simes' line)**.** *Let $\mathcal{T} : F \mapsto \mathcal{U}(F, \mathcal{A}(F))$ a threshold function based on Simes' line. If the level function $\mathcal{A}$ is Hadamard-differentiable at $G$, tangentially to $C[0,1]$, and satisfies existence Condition C.10, then*

$$\sqrt{m}\left(\mathsf{FDP}_m(\mathcal{T}(\widehat{\mathbb{G}}_m)) - \pi_0\mathcal{A}(G)\right) \rightsquigarrow \pi_0\mathcal{A}(G)\left(\frac{\mathbb{Z}_0(\tau^\star)}{\tau^\star} + \frac{\dot{\mathcal{A}}_G(\mathbb{Z})}{\mathcal{A}(G)}\right),$$

*with $\mathbb{Z} = \pi_0\mathbb{Z}_0 + (1-\pi_0)\mathbb{Z}_1$, where $\mathbb{Z}_0$ and $\mathbb{Z}_1$ are independent Gaussian processes such that $\mathbb{Z}_0 \overset{(d)}{=} \mathbb{B}$ and $\mathbb{Z}_1 \overset{(d)}{=} \mathbb{B} \circ G_1$, and $\mathbb{B}$ is a standard Brownian bridge on $[0,1]$.*

We consider the two types of plug-in procedures most widely used and theoretically justified: Sto02-like procedures (Sto02 [21], STS04 [23]), in which $\pi_0$ is estimated by $\frac{1-\widehat{\mathbb{G}}_m(\lambda)}{1-\lambda}$ or a slight variant, and the BKY06 procedure [2], in which an upper bound for $\pi_0$ is derived from a first application of the classical BH95 procedure.

**Procedure** Sto02**.**

**Definition 4.13** (Procedure Sto02 [21])**.** *Procedure* Sto02 *is the multiple testing procedure with threshold function*

$$\mathcal{T}^{\mathsf{Sto02}(\lambda)}(F) = \sup\left\{u \in [0,1], F(u) \geq \frac{u}{\alpha}\frac{1-F(\lambda)}{1-\lambda}\right\} .$$

*The level function of this procedure is therefore*

$$\mathcal{A}(F) = \frac{\alpha}{\overline{\pi_0}^F(\lambda)},$$

*with*

$$\overline{\pi_0}^F(\lambda) = \frac{1-F(\lambda)}{1-\lambda} .$$

$\overline{\pi_0}^G(\lambda)$ *will simply be denoted by* $\overline{\pi_0}(\lambda)$*.*

**Condition C.11** (Condition C.2 for procedure Sto02$(\lambda)$)**.** *The target* FDR *level* $\alpha$ *is greater than* $\overline{\pi_0}(\lambda)\,\alpha^\star$*, where* $\alpha^\star$ *is the critical value of the* BH95 *procedure.*

This procedure is known to provide asymptotic control of FDR at level $\alpha$ [21], but does not necessarily control FDR at level $\alpha$ for finite sample size. This led to the definition of a modification of the Sto02 procedure that does control FDR even for finite sample size [23]:

**Definition 4.14** (Procedure STS04$(\lambda)$ [23])**.** *Procedure* STS04$(\lambda)$ *rejects p-values smaller than*

$$\mathcal{T}_m^{\mathsf{STS04}(\lambda)}(\widehat{\mathbb{G}}_m) = \sup\left\{u \in [0,\lambda], \widehat{\mathbb{G}}_m(u) \geq \frac{u}{\alpha}\frac{1+\frac{1}{m}-\widehat{\mathbb{G}}_m(\lambda)}{1-\lambda}\right\} .$$

According to Proposition 3.6, procedures STS04$(\lambda)$ and Sto02$(\lambda)$ are asymptotically equivalent provided that Conditions C.9 and C.11 hold (see Proposition 7.16 page 1104 for a formal proof).

**Theorem 4.15** (Asymptotic properties of the Sto02/STS04 procedure)**.** *Let* $\lambda \in (0,1)$*, and* $\tau^\star = \mathcal{T}^{\mathsf{Sto02}(\lambda)}(G)$*. Under existence Condition C.11, we have*

$$\sqrt{m}\left(\mathsf{FDP}_m(\mathcal{T}^{\mathsf{Sto02}}(\widehat{\mathbb{G}}_m)) - \frac{\pi_0}{\overline{\pi_0}(\lambda)}\alpha\right) \rightsquigarrow X^{\mathsf{Sto02}},$$

*where*

$$X^{\mathsf{Sto02}} = \frac{\pi_0\alpha}{\overline{\pi_0}(\lambda)}\left(\frac{\mathbb{Z}_0(\tau^\star)}{\tau^\star} + \frac{\mathbb{Z}(\lambda)}{1-G(\lambda)}\right),$$

*with* $\mathbb{Z} = \pi_0\mathbb{Z}_0 + (1-\pi_0)\mathbb{Z}_1$, *where* $\mathbb{Z}_0$ *and* $\mathbb{Z}_1$ *are independent Gaussian processes such that* $\mathbb{Z}_0 \overset{(d)}{=} \mathbb{B}$ *and* $\mathbb{Z}_1 \overset{(d)}{=} \mathbb{B} \circ G_1$, *and* $\mathbb{B}$ *is a standard Brownian bridge on* $[0,1]$. $X^{\mathsf{Sto02}}$ *is therefore a centered Gaussian random variable, with variance*

$$\left(\frac{\pi_0\alpha}{\overline{\pi_0}(\lambda)}\right)^2 \left\{\frac{1-\tau^\star}{\tau^\star} + \frac{\operatorname{Var}\mathbb{Z}(\lambda)}{(1-G(\lambda))^2} + 2\pi_0\frac{\tau^\star \wedge \lambda - \tau^\star\lambda}{\tau^\star(1-G(\lambda))}\right\},$$

*where*

$$\operatorname{Var}\mathbb{Z}(\lambda) = \pi_0^2\lambda(1-\lambda) + (1-\pi_0)^2 G_1(\lambda)(1-G_1(\lambda)).$$

**Corollary 4.16.** *If* $\tau^\star \leq \lambda$,

$$\operatorname{Var} X^{\mathsf{Sto02}}(\lambda) = \left(\frac{\pi_0\alpha}{\overline{\pi_0}(\lambda)}\right)^2 \left\{\frac{1-\tau^\star}{\tau^\star} + \frac{\operatorname{Var}\mathbb{Z}(\lambda)}{(1-G(\lambda))^2} + 2\frac{\pi_0}{\overline{\pi_0}(\lambda)}\right\}.$$

As $\overline{\pi_0}(\lambda) = \pi_0 + (1-\pi_0)\frac{1-G_1(\lambda)}{1-\lambda}$, we have, for any $\lambda \leq \lambda'$, $\pi_0 \leq \overline{\pi_0}(\lambda') \leq \overline{\pi_0}(\lambda) \leq 1$, BH95o $\gg$ Sto02($\lambda'$) $\gg$ Sto02($\lambda$) $\gg$ BH95.

**Procedure BKY06.** Letting $\beta = \frac{\alpha}{1+\alpha}$, procedure BKY06 involves applying procedure BH95 at level $\frac{\beta}{1-R(\beta)/m}$, where $R(\beta)$ is the number of hypotheses rejected by a first application of the BH95 procedure at level $\beta$. We shall consider a recently proposed generalization of this procedure [4], in which procedure BH95 is applied at level $\frac{1-\lambda}{1-R(\lambda)/m}\alpha$. The original BKY06 procedure corresponds to $\lambda = \frac{\alpha}{1+\alpha}$.

**Definition 4.17** (Procedure BKY06($\lambda$)[2])**.** *Let* $\lambda \in [0,1)$, *and*

$$\mathcal{A}(F) = \alpha\frac{1-\lambda}{1-F(\mathcal{U}(F,\lambda))},$$

*where*

$$\mathcal{U}(F,\lambda) = \sup\left\{u \in [0,1], F(u) \geq \frac{u}{\lambda}\right\}.$$

*The threshold function of procedure* BKY06($\lambda$)*is defined for any* $F \in D[0,1]$ *by* $\mathcal{T}^{\mathsf{BKY06}(\lambda)}(F) = \mathcal{U}(F, \mathcal{A}(F))$, *that is,*

$$\mathcal{T}^{\mathsf{BKY06}(\lambda)}(F) = \sup\left\{u \in [0,1], F(u) \geq \frac{u}{\alpha}\,\frac{1-F(\mathcal{U}(F,\lambda))}{1-\lambda}\right\}.$$

*Remark* 4.18. As the proportion $R(\lambda)/m$ of hypotheses rejected by procedure BH95 at level $\lambda$ equals $\widehat{\mathbb{G}}_m(\mathcal{U}(\widehat{\mathbb{G}}_m,\lambda))$ (see Proposition 7.8 page 1098), $\mathcal{A}(\widehat{\mathbb{G}}_m)$ may be written as $\alpha\frac{1-\lambda}{1-R(\lambda)/m}$.

*Remark* 4.19. The exact definition of procedure BKY06($\lambda$) adds $1/m$ to the denominator of the level function:

$$\mathcal{A}(F) = \alpha\frac{1-\lambda}{1+\frac{1}{m}-F(\mathcal{U}(F,\lambda))},$$

which permits proving that this procedure controls FDR for finite sample size [4]. According to Proposition 3.6, these two procedures are asymptotically equivalent, so we will use Definition 4.17.

As procedure BKY06($\lambda$) is based on two successive applications of procedure BH95, at level $\lambda$ and $\alpha(1-\lambda)$, Condition C.2 holds if and only if Condition C.8 holds and $\lambda > \alpha^\star$.

**Condition C.12.** *The parameter $\lambda$ satisfies $\lambda > \alpha^\star$, where $\alpha^\star$ is the critical value of the* BH95 *procedure.*

**Theorem 4.20** (Asymptotic properties of the BKY06($\lambda$) procedure)**.** *Let $\alpha \in [0,1]$, and $\lambda \in [0,1)$. Let $u(\lambda) = \mathcal{U}(G,\lambda)$ be the asymptotic threshold of the* BH95 *procedure applied at level $\lambda$, and $\tau^\star = \mathcal{T}^{\mathsf{BKY06}(\lambda)}(G)$. Under existence Conditions C.8 and C.12,*

$$\sqrt{m}\left(\mathsf{FDP}_m(\mathcal{T}^{\mathsf{BKY06}(\lambda)}(\widehat{\mathbb{G}}_m)) - \frac{\pi_0\alpha(1-\lambda)}{1-G(u(\lambda))}\right) \rightsquigarrow X^{\mathsf{BKY06}(\lambda)},$$

*where*

$$X^{\mathsf{BKY06}(\lambda)} = \frac{\pi_0\alpha(1-\lambda)}{1-G(u(\lambda))}\left(\frac{\mathbb{Z}_0(\tau^\star)}{\tau^\star} + \frac{1}{1-\alpha(1-\lambda)g(u(\lambda))}\frac{\mathbb{Z}(u(\lambda))}{1-G(u(\lambda))}\right),$$

*with $\mathbb{Z} = \pi_0\mathbb{Z}_0 + (1-\pi_0)\mathbb{Z}_1$, where $\mathbb{Z}_0$ and $\mathbb{Z}_1$ are independent Gaussian processes such that $\mathbb{Z}_0 \overset{(d)}{=} \mathbb{B}$ and $\mathbb{Z}_1 \overset{(d)}{=} \mathbb{B}\circ G_1$, and $\mathbb{B}$ is a standard Brownian bridge on $[0,1]$.*

As $u(\lambda)$ is the asymptotic threshold of the BH95 procedure applied at level $\lambda$, we have $G(u(\lambda)) = u(\lambda)/\lambda$, $u(\lambda) \leq \lambda$. Therefore, $\frac{1-\lambda}{1-G(u(\lambda))} \leq \frac{1-\lambda}{1-G(\lambda)}$, and the asymptotic level of procedure BKY06($\lambda$) is less than $\alpha$ because $\frac{1-G(\lambda)}{1-\lambda} \geq \pi_0$.

However, as for procedure BR08($\lambda$), procedure BKY06($\lambda$) need not be more powerful than BH95: a comparison of the asymptotic FDR for these two procedures shows that situations in which BKY06($\lambda$) $\gg$ BH95 are characterized by $G(u_\lambda) \geq \lambda$, that is, $G^2(u_\lambda) \geq u_\lambda$ because $G(u_\lambda) = u_\lambda/\lambda$. Hence, BH95 $\gg$ BKY06($\lambda$) corresponds to situations in which $G$ is too close to the Uniform distribution. For example, if $G(x) \leq \sqrt{x}$ for all $x \in [0,1]$, then for any $\lambda \in [0,1]$, BH95 $\gg$ BKY06($\lambda$).

## 5. Connection between one-stage and two-stage adaptive procedures

We have introduced two types of FDR controlling procedures generalizing the BH95 procedure: two-stage adaptive (plug-in) procedures explicitly incorporate an estimate of $\pi_0$ into the standard BH95 procedure, whereas one-stage adaptive procedures do not explicitly use such an estimate, but still provide tighter FDR control than the BH95 procedure.

We will now investigate connections between one-stage and two-stage adaptive procedures, which naturally appear when using the formalism of threshold functions: with a striking symmetry, the threshold of procedure BR08($\lambda$) may be interpreted as a fixed point of an iterated BKY06($\lambda$) procedure, whereas the threshold of procedure FDR08 may be interpreted as a fixed point of an iterated Sto02($\lambda$) procedure. We provide heuristic reasons for these connections in section 5.1; in section 5.2 we present general results for the connection between one-stage and two-stage adaptive procedures, and derive consequences for the connection between procedures Sto02($\lambda$) and FDR08 on the one hand, and between procedures BKY06($\lambda$) and BR08($\lambda$) on the other hand.

### 5.1. Heuristics

**Procedures** BKY06($\lambda$) **and** BR08($\lambda$). The BKY06 procedure was designed to derive an approximate upper bound for $\pi_0$ from a first application of procedure BH95, and to use this upper bound in a second application of the BH95 procedure, leading to less conservative FDR control. For $\lambda \in [0,1)$, the threshold function of the BKY06($\lambda$) procedure is defined by

$$\mathcal{T}^{\mathsf{BKY06}(\lambda)}(F) = \sup\left\{u \in [0,1], F(u) \geq \frac{u}{\alpha}\frac{1-F(\mathcal{U}(F,\lambda))}{1-\lambda}\right\},$$

where $\mathcal{U}(F,\lambda) = \sup\left\{u \in [0,1], F(u) \geq \frac{u}{\lambda}\right\}$. It therefore seems natural to iterate this process, using the number of rejections at the second application to find a less conservative upper bound for $\pi_0$, and to use this new upper bound in a third application of the BH95 procedure, and so on. Based on this idea, Benjamini *et al.* suggested defining a *multi-stage procedure* for the particular situation in which $\lambda = \frac{\alpha}{1+\alpha}$ [2]. In our framework, this iterative process suggests the introduction of a *fixed-point procedure* defined for any $F \in D[0,1]$ by:

$$\mathcal{T}_\infty^{\mathsf{BKY06}(\lambda)}(F) = \sup\left\{u \in [0,1], F(u) \geq \frac{u}{\alpha}\frac{1-F(u)}{1-\lambda}\right\}.$$

The term *fixed-point procedure* refers to the following property of the corresponding asymptotic threshold $\tau_\infty^\star = \mathcal{T}_\infty^{\mathsf{BKY06}(\lambda)}(G)$. Let us suppose that $\tau_\infty^\star$ is the threshold obtained at a given stage of the abovementioned iteration process. As $G(\tau_\infty^\star) = \tau_\infty^\star(1-G(\tau_\infty^\star))/\alpha(1-\lambda)$, $\tau_\infty^\star$ is also the asymptotic threshold at the next stage, and is thus a fixed point of the iteration process. It turns out that *this fixed-point procedure is the* BR08($\lambda$) *procedure* investigated in section 4.2: $F(u) \geq \frac{u}{\alpha(1-\lambda)}(1-F(u))$ may be written as $F(u) \geq \frac{u}{\alpha(1-\lambda)+u}$, and the right-hand side is the rejection curve $b_\alpha^\lambda$ of the BR08($\lambda$) procedure.

**Procedures** Sto02($\lambda$) **and** FDR08($\lambda$). The same idea may be adapted to procedure Sto02($\lambda$), which is defined for $0 \leq \lambda < 1$ by the threshold function

$$\mathcal{T}^{\mathsf{Sto02}(\lambda)}(F) = \sup\left\{u \in [0,1], F(u) \geq \frac{u}{\alpha}\frac{1-F(\lambda)}{1-\lambda}\right\}.$$

If $\widehat{\tau}_\lambda = \mathcal{T}^{\mathsf{Sto02}(\lambda)}(\widehat{\mathbb{G}}_m)$ denotes the empirical threshold of procedure $\mathsf{Sto02}(\lambda)$, one may use $\widehat{\tau}_\lambda$ to estimate $\pi_0$, that is, calculate the threshold given by procedure $\mathsf{Sto02}(\widehat{\tau}_\lambda)$, and so on. This suggests that an associated *fixed-point procedure* could be defined as

$$\mathcal{T}_\infty^{\mathsf{Sto02}}(F) = \sup\left\{u \in [0,1], F(u) \geq \frac{u}{\alpha}\frac{1-F(u)}{1-u}\right\}.$$

Again, the term *fixed-point procedure* refers to the fact that if $\tau_\infty^\star = \mathcal{T}(G)$ is used as a new $\lambda$ to estimate $\pi_0$ in procedure $\mathsf{Sto02}(\lambda)$, then the asymptotic threshold of procedure $\mathsf{Sto02}(\lambda)$ is also $\tau_\infty^\star$, which is therefore a fixed point of the iteration process. It turns out that *this fixed-point procedure is the* FDR08 *procedure* investigated in section 4.2: $F(u) \geq \frac{u}{\alpha}\frac{1-F(u)}{1-u}$ may be written as $F(u) \geq \frac{u}{\alpha+(1-\alpha)u}$, and the right-hand side is the rejection curve $f_\alpha$ of the FDR08 procedure.

### *5.2. Formal connections*

We present a general result concerning connections between one-stage and two-stage adaptive procedures, providing a formal justification for the connections mentioned in section 5.1, and accounting for their symmetry. This result is based on the following assumption concerning the threshold function of the one-stage adaptive procedure:

**Condition C.13.** *There is a curve* $c_\alpha : D[0,1] \times [0,1]$ *such that the threshold function* $\mathcal{T}$ *may be written as*

$$\mathcal{T}(F) = \sup\left\{u \in [0,\lambda], F(u) \geq c_\alpha(F,u)\right\},$$

*where* $u \mapsto c_\alpha(G,u)/u$ *is non increasing on* $[0,\lambda]$.

*Remark* 5.1. In Condition C.13, $c_\alpha(F,\cdot)$ is *not* the rejection curve of procedure $\mathcal{T}$, because it depends on $F$. For example, for procedure FDR08, we will use

$$c_\alpha(F,u) = \frac{u}{\alpha}\,\frac{1-F(u)}{1-u}.$$

Theorem 5.2 shows that we can associate with a one-stage adaptive procedure fulfilling Condition C.13 a two-stage adaptive procedure with linear rejection curve, and level function given by

$$\mathcal{A}(F) = \frac{t}{c_\alpha(F,t)},$$

for fixed $t \in (0,1)$. The asymptotic threshold of the one-stage procedure may then be interpreted as the fixed point of iterations of the two-stage procedure.

**Theorem 5.2** (Connection between one-stage and two-stage adaptive procedures)**.** *Let* $\lambda \in (0,1)$. *Let us consider a multiple testing procedure with a threshold function* $\mathcal{T}$ *that may be written as*

$$\mathcal{T}(F) = \sup\left\{u \in [0,\lambda], F(u) \geq c_\alpha(F,u)\right\}$$

*for any $F \in D[0,1]$. Let $\mathcal{T}_t$ be the threshold function defined by*

$$\mathcal{T}_t(F) = \sup\left\{u \in [0,1], F(u) \geq \frac{c_\alpha(F,t)}{t} u\right\},$$

*for any $t \in (0,1)$ and any $F \in D[0,1]$. Let us assume that existence Condition C.2 and uniqueness Condition C.3 hold for procedure $\mathcal{T}$, and that, for any $t \in (0,1)$, existence Condition C.2 holds for procedure $\mathcal{T}_t$. Let $\tau^\star = \mathcal{T}(G)$ and $\tau(t) = \mathcal{T}_t(G)$ be the asymptotic thresholds of procedures $\mathcal{T}$ and $\mathcal{T}_t$, respectively. If $c_\alpha$ satisfies Condition C.13, we have*

*(i)* *for any $t \in (0,\lambda]$,*

$$\begin{cases} t \leq \tau^\star \Rightarrow \tau(t) \in [t, \tau^\star] \\ t \geq \tau^\star \Rightarrow \tau(t) \in [\tau^\star, t] \end{cases}.$$

*(ii)* *Let $t \in (0,\lambda]$. Define the sequence $(t_n) \in [0,1]^{\mathbb{N}}$ by $t_0 = t$, and $t_{i+1} = \tau(t_i)$ for $i \in \mathbb{N}$. Then*

$$\lim_{n\to\infty} t_n = \tau^\star.$$

**Corollary 5.3** (Asymptotic power comparison)**.** *With the same notation and under the same conditions, the following assertions are equivalent:*

*(i)* *Procedure $\mathcal{T}_t$ is asymptotically more powerful than procedure $\mathcal{T}$*
*(ii)* $\tau(t) > \tau^\star$
*(iii)* $t > \tau^\star$
*(iv)* $t > \tau(t)$

In the remainder of this section, we use Theorem 5.2 to characterize the connection between the abovementioned procedures.

**Procedures Sto02($\lambda$) and FDR08($\lambda$).** Theorem 5.4 gives the convergence of the process consisting of the recursive use of the asymptotic threshold of procedure Sto02($\lambda$) as a new $\lambda$. It holds under the same regularity conditions as those required to obtain the asymptotic distribution of procedure FDR08.

**Theorem 5.4** (Connection between procedures Sto02($\lambda$) and FDR08)**.** *Let $\kappa = \frac{\alpha(1-\pi_0)}{\pi_0(1-\alpha)}$, and*

$$\tau^\star = \sup\left\{u \in [0,\kappa], G(u) \geq \frac{u}{\alpha}\frac{1-G(u)}{1-u}\right\}$$

*be the asymptotic threshold of the* FDR08 *procedure. For $u \in [0,1]$, let*

$$\tau(u) = \sup\left\{u \in [0,1], G(u) \geq \frac{u}{\alpha}\frac{1-G(\lambda)}{1-\lambda}\right\}$$

*be the asymptotic threshold of procedure* Sto02($u$)*. For any $t \in (0,1)$, define the sequence $(t_n) \in [0,1]^{\mathbb{N}}$ by $t_0 = t$, and $t_{i+1} = \tau(t_i)$ for $i \in \mathbb{N}$. Let us assume that uniqueness Condition C.3 holds for procedure* FDR08*, and that the target* FDR *level $\alpha$ satisfies existence Conditions C.4 and C.7. Then,*

$$\lim_{n\to\infty} t_n = \tau^\star .$$

**Corollary 5.5** (Asymptotic power comparison — Sto02($\lambda$) vs FDR08)**.** *With the same notation and under the same conditions, procedure* Sto02($\lambda$) *is asymptotically more powerful than procedure* FDR08 *if and only if $\lambda > \tau(\lambda)$.*

When using procedure Sto02($\lambda$) in practice, we would not want any of the rejected hypotheses to be incorporated into the estimation of $\pi_0$. Thus, the empirical rejection threshold $\mathcal{T}^{\mathsf{Sto02}(\lambda)}(\widehat{\mathbb{G}}_m)$ should be less than $\lambda$. In such situations, as $\mathcal{T}^{\mathsf{Sto02}(\lambda)}(\widehat{\mathbb{G}}_m)$ converges at rate $1/\sqrt{m}$ to $\tau(\lambda) = \mathcal{T}^{\mathsf{Sto02}(\lambda)}(G)$, procedure Sto02($\lambda$) is probably more powerful than procedure FDR08 according to Corollary 5.5.

**Procedures BKY06($\lambda$) and BR08($\lambda$).** Theorem 5.6 characterizes the connection between procedure BKY06($\lambda$) and procedure BR08($\lambda$). It holds under the same regularity conditions as those required to obtain the asymptotic distribution of procedure BR08($\lambda$).

Let $\tau^\star = \mathcal{T}^{\mathsf{BR08}(\lambda)}(G)$ be the asymptotic threshold of the BR08($\lambda$) procedure. Under uniqueness Condition C.3 and existence Conditions C.8, C.9 and C.12, $(t_n)$ is non decreasing, and converges to $\tau^\star$.

**Theorem 5.6** (Connection between procedures BKY06($\lambda$) and BR08($\lambda$))**.** *Let $\lambda \in (0,1)$. For $F \in D[0,1]$ and $\beta \in [0,1]$, let*

$$\mathcal{U}(F,\beta) = \sup\left\{u \in [0,1], F(u) \geq \frac{u}{\beta}\right\} .$$

*For any $u \in [0,1)$, let $\tau(u) = \mathcal{U}\left(G, \frac{\alpha(1-\lambda)}{1-G(u)}\right)$. With this notation, $\mathcal{T}(G) = \tau(u(\lambda))$ is the asymptotic threshold of the* BKY06($\lambda$) *procedure, where $u(\lambda) = \mathcal{U}(G,\lambda)$. Let*

$$\tau^\star = \sup\left\{u \in [0,\lambda], G(u) \geq \frac{u}{\alpha(1-\lambda)+u}\right\}$$

*be the asymptotic threshold of the* BR08($\lambda$) *procedure. Define the sequence $(t_n) \in [0,1]^{\mathbb{N}}$ by $t_0 = u(\lambda)$, and $t_{i+1} = \tau(t_i)$ for $i \in \mathbb{N}$. Let us assume that uniqueness Condition C.3 holds for procedure* BR08($\lambda$)*, and that the target* FDR *level $\alpha$ satisfies existence Conditions C.8 and C.9. Then*

$$\lim_{n\to\infty} t_n = \tau^\star .$$

**Corollary 5.7** (Asymptotic power comparison — BKY06($\lambda$) vs BR08($\lambda$))**.** *With the same notation and under the same conditions, procedure* BR08($\lambda$) *is asymptotically more powerful than procedure* BKY06($\lambda$) *if and only if the asymptotic threshold $\tau^\star$ of procedure* BR08($\lambda$) *satisfies $\tau^\star \geq \lambda - \alpha(1-\lambda)$.*

For example, setting $\lambda$ to a value less than $\frac{\alpha}{1+\alpha}$, corresponding to the original BKY06 procedure [2], ensures that the associated BR08($\lambda$) procedure is asymptotically more powerful than the associated BKY06($\lambda$) procedure.

## 6. Concluding remarks

This paper demonstrates the power and flexibility of the formalism of threshold functions, making it possible to derive the asymptotic properties of well known FDR controlling procedures with their associated regularity conditions, and to identify and characterize novel connections between one-stage and two-stage adaptive procedures. These results are summarized in Table 1. We should recall that the threshold function associated with the level function $\mathcal{A}$ and rejection curve $r_\alpha = r(\alpha, \cdot)$ is defined by

$$\mathcal{T}(F) = \sup \{u \in [0,1], F(u) \geq r(\mathcal{A}(F), u)\} .$$

By definition, the level function $\mathcal{A}$ equals $\alpha$ for one-stage procedures, and the rejection curve of Simes' line-based procedures is $r_\alpha : u \mapsto u/\alpha$.

**Regularity conditions.** For one-stage adaptive procedures FDR08 and BR08($\lambda$), the uniqueness Condition C.3 has to be assumed (cf. Table 1): as the rejection curve is not linear, the interior right crossing point is not necessarily unique; in practice the uniqueness condition holds except in pathological situations. For Simes' line-based procedures BH95, Sto02($\lambda$) and BKY06($\lambda$), existence Condition C.2 holds provided that the slope of the distribution function exceeds a certain threshold at the origin (that is, that there is no criticality phenomenon). For one-stage adaptive procedures, it is also required that the rejection curve $r_\alpha$ ends below the distribution function $G$, which corresponds to Condition C.7 for procedure FDR08, and Condition C.9 for procedure BR08($\lambda$).

The criticality phenomenon studied by [5] is intrinsic to the multiple testing problem, and not specific to a given procedure, as the minimum attainable pFDR level $\beta^\star = \inf_{t>0} \mathsf{pFDR}(t)$ depends solely on the parameters of the model [5]. When $\beta^\star = 0$, say for the Gaussian location problem, *there is no criticality phenomenon for any procedure*: $\alpha^\star = 0$, and all existence Conditions concerning the behavior of the distribution function $G$ close to 0 are fulfilled for any procedure, and for any target FDR level $\alpha$. When $\beta^\star > 0$, say for the Laplace

TABLE 1
*Comparison of* FDR *controlling procedures, characterized by their level function $\mathcal{A}$ and their rejection curve $r_\alpha$. Conditions for the existence and uniqueness of an interior right crossing point are recalled, together with the corresponding* pFDR *relative to that of the* BH95 *procedure: $\pi_0\alpha$*

| Name | BH95 [1] | FDR08 [7] | BR08($\lambda$) [4] | Sto02($\lambda$) [21] | BKY06($\lambda$) [2] |
|---|---|---|---|---|---|
| $\mathcal{A}(F)/\alpha$ | 1 | 1 | 1 | $\frac{1-\lambda}{1-F(\lambda)}$ | $\frac{1-\lambda}{1-\widehat{\mathbb{G}}_m(u_\lambda)}$ $^{(a)}$ |
| $r_\alpha(u)$ | $u/\alpha$ | $\frac{u}{\alpha+(1-\alpha)u}$ $^{(b)}$ | $\frac{u}{\alpha(1-\lambda)+u}$ $^{(c)}$ | $u/\alpha$ | $u/\alpha$ |
| Existence | C.4 | C.4 & C.7 $^{(d)}$ | C.8 & C.9 $^{(d)}$ | C.11 | C.12 |
| Uniqueness | — | C.3 | C.3 | — | — |
| pFDR$/\pi_0\alpha$ | 1 | $\frac{1-\tau^\star_{\mathsf{FDR08}}}{1-G\left(\tau^\star_{\mathsf{FDR08}}\right)}$ | $\frac{1-\lambda}{1-G\left(\tau^\star_{\mathsf{BR08}}\right)}$ | $\frac{1-\lambda}{1-G(\lambda)}$ | $\frac{1-\lambda}{1-G(u_\lambda)}$ $^{(a)}$ |

$(a)$ : $u_\lambda$ is the asymptotic threshold of the BH95 procedure at target level $\lambda$; $(b)$ : truncated at $\frac{\alpha(1-\pi_0)}{\pi_0(1-\alpha)}$; $(c)$ : truncated at $\lambda$; $(d)$ : Sufficient (not necessary) conditions.

location problem (Figure 2, page 1076), *there is a criticality phenomenon for every procedure*; however the critical value, that is, the minimum target FDR level for which existence Condition C.2 holds, may depend on the procedure, as illustrated by the existence conditions in Table 1.

**Power comparisons.** All procedures are asymptotically conservative, and therefore yield asymptotic FDR below the target level. Procedures FDR08 and Sto02 (and thus STS04) are always more powerful than procedure BH95, but this is not the always the case for procedures BR08($\lambda$)(section 4.2) and BKY06($\lambda$) (section 4.3).

For one-stage adaptive procedures, for any $\lambda \in (0, 1)$ such that the regularity conditions for procedures FDR08 and BR08($\lambda$) hold, FDR08 is asymptotically more powerful than BR08($\lambda$). Indeed, Condition C.9 ensures that the asymptotic thresholds of both procedures are less than $\lambda$. As the rejection curve $f_\alpha$ of procedure FDR08 is smaller than the rejection curve $b_\alpha^\lambda$ of BR08 on $[0, \lambda]$, the asymptotic threshold of procedure FDR08 is greater than that of procedure BR08($\lambda$). However, it should be noted that procedure BR08($\lambda$) does control FDR for a finite number of tested hypotheses, whereas procedure FDR08 does not.

For two-stage adaptive procedures, for any $\lambda \in (0, 1)$ such that the regularity conditions for procedures Sto02($\lambda$) and BKY06($\lambda$) hold, Sto02($\lambda$) (and thus STS04) is asymptotically more powerful than BKY06($\lambda$), as demonstrated by the corresponding asymptotic FDR levels in Table 1: as $u_\lambda \leq \lambda$, we have $\frac{1-\lambda}{1-G(u_\lambda)} \leq \frac{1-\lambda}{1-G(\lambda)}$. This suggests that procedure STS04($\lambda$) is preferable to procedure BKY06($\lambda$) in practice. This recommendation should be balanced against the choice of $\lambda$ and the desired robustness to dependence between null hypotheses. Based on a simulation study, procedure Sto02($\alpha$) was recently reported to be much more robust to positive dependence between null hypotheses than procedure Sto02(1/2) [4], which is still a standard choice in practical implementations, such as the SAM (Significance Analysis of Microarrays) software [24].

**Towards optimality.** This comparison raises the question of whether the formalism of threshold functions can be used to derive procedures more powerful than those studied here. One possible approach consists of trying to improve the estimation of $\pi_0$ to build a procedure closer to the Oracle BH95 procedure, as discussed in [11]. However, consistent estimators of $\pi_0$ have slower convergence rates than $1/\sqrt{m}$, resulting in slower convergence rates than $1/\sqrt{m}$ for the associated FDP. This may be illustrated by the influence of $\lambda$ on procedure Sto02($\lambda$): the larger $\lambda$, the smaller the bias $\mathbb{E}[\widehat{\pi_0}(\lambda)] - \pi_0$, and the larger the variance of $\widehat{\pi_0}(\lambda)$. The question of how to choose $\lambda$ as a function of the number of hypotheses tested and the assumed regularity of $G$ is discussed in another work [14].

Another possibility would be to consider procedures more general than those used in this paper: the BH95o procedure has been shown to give the lowest false non discovery rates (FNR) of the threshold procedures controlling FDR at level $\alpha$ [10]. The question of optimality in a broader family of testing procedures has recently been raised [25]: $Z$ score-based threshold procedures may outperform $p$ value-based threshold procedures, as they make it possible to choose different significance thresholds for positive and negative significance cutoffs. This suggests to extend our framework to $Z$ score-based procedures.

**Confidence intervals.** An interesting practical application of this work concerns the derivation of asymptotic confidence intervals for the FDP of a given procedure. Our results give explicit asymptotic distributions for the attained FDP, but this issue is not straightforward because these distributions depend on unknown quantities, including the proportion $\pi_0$, the asymptotic threshold FDR $\tau^\star$, or the distribution function $G$ and its associated density $g$. These quantities should, in turn, be estimated. Bootstrapping techniques could be used for this purpose; we leave this question for further research.

**Extension to other dependence settings.** We have derived the asymptotic properties of several multiple testing procedures and the associated regularity conditions in the situation in which $p$-values are independent. However, our formalism makes it possible to deal with any dependence situation for which the vector $(\widehat{\mathbb{G}}_{0,m}, \widehat{\mathbb{G}}_{1,m})$ of empirical distribution functions of the $p$-values under the null and alternative hypotheses satisfies Donsker's invariance principle. For example, the form of the asymptotic distributions of the threshold $\mathcal{T}(\widehat{\mathbb{G}}_m)$ and the associated FDP would remain the same in the conditional dependence model recently proposed by Wu [28].

## 7. Proof of main results

### 7.1. Asymptotic FDP: general threshold functions

In this section, we provide proofs for the results of section 3.

#### 7.1.1. Proof of Theorem 3.2

The following lemma will be used in several subsequent proofs.

**Lemma 7.1.** *Let $H \in C[0,1]$, and $H_t$ be a family of functions of $D[0,1]$ that converges to $H$ on $(D[0,1], \|.\|_\infty)$ as $t \to 0$. For any sequence $(u_t)_{t>0}$ of $[0,1]$ that converges to $u \in [0,1]$ as $t \to 0$, we have*

$$\begin{aligned} \lim_{t\to 0} H_t(u_t) &= H(u) \\ \lim_{t\to 0} H_t(u_t^-) &= H(u), \end{aligned}$$

*where $f(x_0^-)$ denotes $\lim_{x\to x_0, x\leq x_0} f(x)$.*

*Proof of Lemma 7.1.* We have

$$\begin{aligned} |H_t(u_t) - H(u)| &\leq |H_t(u_t) - H(u_t)| + |H(u_t) - H(u)| \\ &\leq \|H_t - H\|_\infty + |H(u_t) - H(u)| \end{aligned}$$

and

$$\begin{aligned} |H_t(u_t^-) - H(u)| &\leq |H_t(u_t^-) - H(u_t^-)| + |H(u_t^-) - H(u)| \\ &\leq \|H_t - H\|_\infty + |H(u_t) - H(u)| \end{aligned}$$

as $H$ is continuous. The first term goes to 0 as $t \to 0$ by the convergence of $H_t$ to $H$ on $D[0,1]$, and the second term also tends to 0 by the continuity of $H$, because $\lim_{t\to 0} u_t = u$. □

**Proposition 7.2** (Hadamard differentiability of $\mathcal{V}$ and $\mathcal{R}$)**.** *Under Condition C.1,*

*(i)* *$\mathcal{V}$ is Hadamard-differentiable at $(G_0, G_1)$, tangentially to $C[0,1]^2$, with derivative*

$$\dot{\mathcal{V}}_{G_0,G_1} : (H_0, H_1) \mapsto \pi_0 \dot{\mathcal{T}}_G \left(\pi_0 H_0 + (1-\pi_0) H_1\right) + \pi_0 H_0(\mathcal{T}(G))$$

*(ii)* *$\mathcal{R}$ is Hadamard-differentiable at $G$, tangentially to $C[0,1]$, with derivative*

$$\dot{\mathcal{R}}_G : H \mapsto H(\tau^\star) + g(\tau^\star)\dot{\mathcal{T}}_G(H)$$

*Proof of Proposition 7.2.* (i) Let $(H_0, H_1) \in C[0,1]^2$, and $(H_{0,t}, H_{1,t})_{t>0}$ be a family of functions of $D[0,1]^2$ that converges to $((H_0,H_1), \|.\|_\infty)$ as $t \to 0$. Let $H = \pi_0 H_0 + (1-\pi_0)H_1$, and $H_t = \pi_0 H_{0,t} + (1-\pi_0)H_{1,t}$. We have

$$\mathcal{V}(G_0 + tH_{0,t}, G_1 + tH_{1,t}) - \mathcal{V}(G_0, G_1) = \pi_0(\tau_t^\star - \tau^\star) + \pi_0 t H_{0,t}(\tau_t^\star)$$

where $\tau^\star = \mathcal{T}(G)$ and $\tau_t^\star$ denotes $\mathcal{T}(G + tH_t)$. By the Hadamard differentiability of $\mathcal{T}$ at $G$ tangentially to $C[0,1]$, we have, as $H = \pi_0 H_0 + (1 - \pi_0)H_1$ is continuous at $\tau^\star$,

$$\tau_t^\star - \tau^\star = t\left(\dot{\mathcal{T}}_G(H) + o(1)\right)$$

In order to conclude, we notice that

$$\lim_{t\to 0} H_{0,t}(\tau_t^\star) \to H_0(\tau^\star)$$

according to Lemma 7.1, which concludes the proof.

(ii) Let $H \in C[0,1]$, and $H_t$ be a family of functions of $D[0,1]$ that converges to $H$ on $(D[0,1], \|.\|_\infty)$ as $t \to 0$. We have

$$\begin{aligned} \mathcal{R}(G + tH_t) &= (G + tH_t)\mathcal{T}(G + tH_t) \\ &= G(\mathcal{T}(G + tH_t)) + tH_t(\mathcal{T}(G + tH_t)) \end{aligned}$$

By the Hadamard differentiability of $\mathcal{T}$ at $G$ tangentially to $C[0,1]$, we have

$$\mathcal{T}(G + tH_t) = \mathcal{T}(G) + t\left(\dot{\mathcal{T}}_G(H) + o(1)\right)$$

so that applying Taylor's formula to $G$ at $\mathcal{T}(G)$ yields

$$G(\mathcal{T}(G + tH_t)) = G(\mathcal{T}(G)) + t\left(\dot{\mathcal{T}}_G(H) + o(1)\right) g(\mathcal{T}(G)) + o(t)\,.$$

For the second term, Lemma 7.1 ensures that

$$\lim_{t\to 0} H_t(\mathcal{T}(G+tH_t)) = H(\mathcal{T}(G))$$

because $\mathcal{T}(G+tH_t)$ converges to $\mathcal{T}(G)$ and $H_t$ converges to $H$ on $(D[0,1], \|.\|_\infty)$. Finally, we have

$$\lim_{t\to 0} \frac{\mathcal{R}(G+tH_t)-\mathcal{R}(G)}{t} = H(\tau^\star) + g(\tau^\star)\dot{\mathcal{T}}_G(H)$$

because $\tau^\star = \mathcal{T}(G)$, which concludes the proof. □

**Theorem 7.3** (Asymptotic distribution of $(\hat\tau, \hat\nu, \hat\rho)$)**.** *Under Condition C.1,*

$$\sqrt{m}\left(\begin{pmatrix}\widehat{\tau}\\ \widehat{\nu}\\ \widehat{\rho}\end{pmatrix} - \begin{pmatrix}\tau^\star\\ \pi_0\tau^\star\\ r_\alpha(\tau^\star)\end{pmatrix}\right) \rightsquigarrow X\,,$$

*where*

$$X = \begin{pmatrix}1\\ \pi_0\\ g(\tau^\star)\end{pmatrix}\dot{\mathcal{T}}_G(\mathbb{Z}) + \pi_0\begin{pmatrix}0\\1\\1\end{pmatrix}\mathbb{Z}_0(\tau^\star) + (1-\pi_0)\begin{pmatrix}0\\0\\1\end{pmatrix}\mathbb{Z}_1(\tau^\star)\,,$$

*with $\mathbb{Z} = \pi_0\mathbb{Z}_0 + (1-\pi_0)\mathbb{Z}_1$, where $\mathbb{Z}_0$ and $\mathbb{Z}_1$ are independent Gaussian processes such that $\mathbb{Z}_0 \overset{(d)}{=} \mathbb{B}$ and $\mathbb{Z}_1 \overset{(d)}{=} \mathbb{B}\circ G_1$, and $\mathbb{B}$ is a standard Brownian bridge on* $[0,1]$.

*Proof of Theorem 7.3.* We note that

$$\begin{pmatrix}\widehat{\tau}\\ \widehat{\nu}\\ \widehat{\rho}\end{pmatrix} = \Psi(\widehat{\mathbb{G}}_{0,m}, \widehat{\mathbb{G}}_{1,m})$$

where $\Psi : D[0,1]^2 \to \mathbb{R}^3$ is the map defined by

$$\Psi(F_0,F_1) = \begin{pmatrix}\mathcal{T}(\pi_0F_0 + (1-\pi_0)F_1)\\ \mathcal{V}(F_0,F_1)\\ \mathcal{R}(\pi_0F_0+(1-\pi_0)F_1)\end{pmatrix}.$$

We have

$$\Psi(G_0,G_1) = \begin{pmatrix}\tau^\star\\ \pi_0\tau^\star\\ G(\tau^\star)\end{pmatrix}.$$

By the Hadamard differentiability of $\mathcal{T}$ at $G = \pi_0G_0 + (1-\pi_0)G_1$ and that of $\mathcal{V}$ at $(G_0,G_1)$, $\Psi$ is Hadamard-differentiable at $(G_0,G_1)$ tangentially to $C[0,1]^2$, with derivative

$$\dot\Psi_{G_0,G_1}(H_0,H_1) = \begin{pmatrix}\dot{\mathcal{T}}_G(H)\\ \dot{\mathcal{V}}_{G_0,G_1}(H_0,H_1)\\ \dot{\mathcal{R}}_G(H)\end{pmatrix}$$

where $H$ denotes $\pi_0 H_0 + (1-\pi_0)H_1$. Therefore Theorem 3.1 yields

$$\sqrt{m}(\Psi(\widehat{\mathbb{G}}_{0,m}, \widehat{\mathbb{G}}_{1,m}) - \Psi(G_0, G_1)) \rightsquigarrow \dot{\Psi}_{G_0,G_1}(\mathbb{Z}_0, \mathbb{Z}_1),$$

According to Proposition 7.2, we have

$$\begin{aligned} \dot{\mathcal{V}}_{G_0,G_1}(\mathbb{Z}_0, \mathbb{Z}_1) &= \pi_0 \dot{\mathcal{T}}_G(\mathbb{Z}) + \pi_0 \mathbb{Z}_0(\tau^\star) \\ \dot{\mathcal{R}}_G(\mathbb{Z}) &= g(\tau^\star)\dot{\mathcal{T}}_G(\mathbb{Z}) + \mathbb{Z}(\tau^\star) \end{aligned}$$

with $\mathbb{Z} = \pi_0 \mathbb{Z} + (1-\pi_0)\mathbb{Z}_1$, so that

$$\begin{aligned} X &= \dot{\Psi}_{G_0,G_1}(\mathbb{Z}_0, \mathbb{Z}_1) \\ &= \begin{pmatrix} 1 \\ \pi_0 \\ g(\tau^\star) \end{pmatrix} \dot{\mathcal{T}}_G(\mathbb{Z}) + \pi_0 \begin{pmatrix} 0 \\ 1 \\ 1 \end{pmatrix} \mathbb{Z}_0(\tau^\star) + (1-\pi_0) \begin{pmatrix} 0 \\ 0 \\ 1 \end{pmatrix} \mathbb{Z}_1(\tau^\star) \end{aligned}$$

, which concludes the proof. □

*Proof of Theorem 3.2.* $(i)$ is a direct consequence of Theorem 7.3:

$$\sqrt{m}\left(\mathcal{T}(\widehat{\mathbb{G}}_m) - \tau^\star)\right) \rightsquigarrow \dot{\mathcal{T}}_G(\mathbb{Z}),$$

For $(ii)$ and $(iii)$, we note that as $\tau^\star > 0$ (by Condition C.1), $\widehat{\tau} = \mathcal{T}(\widehat{\mathbb{G}}_m)$ is bounded away from 0 if $m$ is sufficiently large, with probability 1. Specifically, there exist $t_0 \in (0,1)$ and $m_0 \in \mathbb{N}$ such that

$$\mathbb{P}\left(\forall m \geq m_0, \widehat{\tau} > t_0\right) = 1.$$

Therefore, as $\widehat{\mathbb{G}}_m$ is non decreasing, and as $\widehat{\mathbb{G}}_m(t_0) \to G(t_0) > 0$, the proportion $\widehat{\rho} = \widehat{\mathbb{G}}_m(\widehat{\tau})$ of rejections by procedure $\mathcal{T}$ is bounded away from 0 with probability 1. Thus, $\mathbb{P}\left(R_m(\widehat{\tau}) > 0\right) = 1$, where $R_m(\widehat{\tau}) = m\widehat{\mathbb{G}}_m(\widehat{\tau})$ is the number of rejections at threshold $\widehat{\tau}$.

As a first consequence, as $\mathsf{FDR}(t) = \mathsf{p}(t)\mathbb{P}\left(R(t) > 0\right)$, we have $\mathsf{FDR}_m(\widehat{\tau}) = \mathsf{p}(\widehat{\tau})$ for a sufficiently large $m$, which proves $(ii)$ as $\mathsf{p}(\widehat{\tau})$ converges almost surely to $\mathsf{p}(\tau^\star)$.

As a second consequence, letting $\gamma : \mathbb{R}_+ \times \mathbb{R}_+^* \to \mathbb{R}$ be defined by $\gamma(x,y) = \frac{x}{y}$, we have $\mathsf{FDP}_m(\mathcal{T}(\widehat{\mathbb{G}}_m)) = \gamma\left(\widehat{\nu}, \widehat{\rho}\right) \mathbf{1}_{\widehat{\rho}>0} = \gamma\left(\widehat{\nu}, \widehat{\rho}\right)$ for a sufficiently large $m$, with probability 1. $\gamma$ is differentiable on $\mathbb{R}_+ \times \mathbb{R}_+^*$, with derivative

$$\dot{\gamma}_{x,y} = \left(\frac{1}{y}, -\frac{x}{y^2}\right).$$

In particular, $\dot{\gamma}_{\pi_0\tau^\star, G(\tau^\star)}(h,k) = \frac{1}{G(\tau^\star)}\left(h - \frac{\pi_0\tau^\star}{G(\tau^\star)}k\right)$. We can therefore derive the asymptotic distribution of $\mathsf{FDP}_m$ from Theorem 7.3 combined with the Delta method [27]. According to Theorem 7.3 we have

$$\sqrt{m}\left(\begin{pmatrix} \widehat{\nu} \\ \widehat{\rho} \end{pmatrix} - \begin{pmatrix} \pi_0\tau^\star \\ G(\tau^\star) \end{pmatrix}\right) \rightsquigarrow \begin{pmatrix} \pi_0\mathbb{Z}_0(\tau^\star) + \pi_0\dot{\mathcal{T}}_G(\mathbb{Z}) \\ \mathbb{Z}(\tau^\star) + g(\tau^\star)\dot{\mathcal{T}}_G(\mathbb{Z}) \end{pmatrix}.$$

Hence, as $\mathsf{FDP}_m(\mathcal{T}(\widehat{\mathbb{G}}_m)) = \gamma(\widehat{\nu}, \widehat{\rho})$ (almost surely), and $\gamma(\pi_0\tau^\star, \tau^\star/\alpha) = \pi_0\alpha$, the Delta method [27] yields

$$\sqrt{m}\left(\mathsf{FDP}_m(\mathcal{T}(\widehat{\mathbb{G}}_m)) - \frac{\pi_0\tau^\star}{G(\tau^\star)}\right) \rightsquigarrow X,$$

where

$$\begin{aligned} X &= \frac{1}{G(\tau^\star)}\left(\pi_0(\mathbb{Z}_0(\tau^\star) + \dot{\mathcal{T}}_G(\mathbb{Z})) - \frac{\pi_0\tau^\star}{G(\tau^\star)}(\mathbb{Z}(\tau^\star) + g(\tau^\star)\dot{\mathcal{T}}_G(\mathbb{Z}))\right) \\ &= \frac{\pi_0\tau^\star}{G(\tau^\star)}\left(\frac{\mathbb{Z}_0(\tau^\star)}{\tau^\star} - \frac{\mathbb{Z}(\tau^\star)}{G(\tau^\star)}\right) + \frac{\pi_0}{G(\tau^\star)}\left(1 - \frac{\tau^\star g(\tau^\star)}{G(\tau^\star)}\right)\dot{\mathcal{T}}_G(\mathbb{Z}) \end{aligned}$$

As $\mathbb{Z} = \pi_0\mathbb{Z}_0 + (1-\pi_0)\,\mathbb{Z}_1$ and $G = \pi_0 G_0 + (1-\pi_0)G_1$, and $\mathsf{p}^\star = \frac{\pi_0 G_0(\tau^\star)}{G(\tau^\star)}$, we have

$$\frac{\mathbb{Z}(\tau^\star)}{G(\tau^\star)} = \mathsf{p}^\star\frac{\mathbb{Z}_0(\tau^\star)}{\tau^\star} + (1-\mathsf{p}^\star)\frac{\mathbb{Z}_1(\tau^\star)}{G_1(\tau^\star)}\,,$$

so that

$$\frac{\mathbb{Z}_0(\tau^\star)}{\tau^\star} - \frac{\mathbb{Z}(\tau^\star)}{G(\tau^\star)} = (1-\mathsf{p}^\star)\left(\frac{\mathbb{Z}_0(\tau^\star)}{\tau^\star} - \frac{\mathbb{Z}_1(\tau^\star)}{G_1(\tau^\star)}\right)\,,$$

which concludes the proof because $\mathsf{p}(t) = \frac{\pi_0 t}{G(t)}$ and $\dot{\mathsf{p}}(t) = \frac{\pi_0}{G(t)}\left(1 - \frac{tg(t)}{G(t)}\right)$. □

#### 7.1.2. *Proof of Proposition 3.6*

Lemma 7.4 states the asymptotic equivalence between a multiple testing procedure defined as a threshold function and a slight modification of this procedure.

**Lemma 7.4.** *Let $\mathcal{T}$ be a threshold function, $\varepsilon = (\varepsilon_m)_{m\in\mathbb{N}}$ and $\mathcal{T}_m^{\varepsilon,H} : D[0,1] \to [0,1]$, such that*

$$\forall F \in D[0,1], \mathcal{T}_m^{\varepsilon,H}(F) = \mathcal{T}(F + \varepsilon_m H)\ .$$

*Let $\mathcal{M}$ be the multiple testing procedure naturally associated with the sequence of thresholds $\mathcal{T}_m^H(\widehat{\mathbb{G}}_m)$. If Condition C.1 holds for $\mathcal{T}$, and if $\varepsilon_m = o\left(\frac{1}{\sqrt{m}}\right)$, then $\mathcal{M}$ is asymptotically equivalent to $\mathcal{T}$ as $m \to +\infty$.*

*Proof of Lemma 7.4.* The proof is based on the idea that, as $\varepsilon_m = o\left(\frac{1}{\sqrt{m}}\right)$, and $\widehat{\mathbb{G}}_m$ converges at rate $\frac{1}{\sqrt{m}}$ to $G$, a modification of $\mathcal{T}$ of the order of $\varepsilon_m$ does not change the asymptotic distribution of the associated FDP, because $\mathcal{T}$ is Hadamard-differentiable. For the sake of simplicity in notation, we prove only that

$$\sqrt{m}\left(\mathcal{T}^{\varepsilon,H}(\widehat{\mathbb{G}}_m) - \mathcal{T}(\widehat{\mathbb{G}}_m)\right) \xrightarrow{P} 0\,.$$

Indeed, as the associated FDP is a Hadamard-differentiable function of the empirical distribution functions under the null and alternative hypotheses $\widehat{\mathbb{G}}_{0,m}$ and $\widehat{\mathbb{G}}_{1,m}$, the arguments developed below can be transposed (but with much more cumbersome notation) to prove that

$$\sqrt{m}\left(\mathsf{FDP}_m(\mathcal{T}^{\varepsilon,H}(\widehat{\mathbb{G}}_m)) - \mathsf{FDP}_m(\mathcal{T}(\widehat{\mathbb{G}}_m))\right) \xrightarrow{P} 0\,.$$

Let $\mathbb{Z}_m = \sqrt{m}\big(\widehat{\mathbb{G}}_m - G\big)$. According to Donsker's invariance principle (Theorem 3.1), $\mathbb{Z}_m$ converges in distribution on $[0,1]$ to a Gaussian process with continuous sample paths. For $Z \in D[0,1]$, let

$$\phi_m(Z) = \sqrt{m}\left(\mathcal{T}^{\varepsilon,H}\left(G + \frac{1}{\sqrt{m}}Z\right) - \mathcal{T}\left(G + \frac{1}{\sqrt{m}}Z\right)\right)\,.$$

We have

$$\phi_m(Z) = \sqrt{m}\left(\mathcal{T}^{\varepsilon,H}\left(G + \frac{1}{\sqrt{m}}Z\right) - \mathcal{T}(G)\right) - \sqrt{m}\left(\mathcal{T}\left(G + \frac{1}{\sqrt{m}}Z\right) - \mathcal{T}(G)\right)\,.$$

According to Condition C.1, $\mathcal{T}$ is Hadamard-differentiable at $G$ tangentially to $C[0,1]$. Therefore, for any sequence $Z_m$ of $D[0,1]$ that converges to $Z \in C[0,1]$, $\sqrt{m}\left(\mathcal{T}\left(G + \frac{1}{\sqrt{m}}Z_m\right) - \mathcal{T}(G)\right)$ converges to $\dot{\mathcal{T}}_G(Z)$. As $\varepsilon_m = o\left(\frac{1}{\sqrt{m}}\right)$, $\sqrt{m}\left(\mathcal{T}^{\varepsilon,H}\left(G + \frac{1}{\sqrt{m}}Z_m\right) - \mathcal{T}(G)\right)$ also converges to $\dot{\mathcal{T}}_G(Z)$.

Thus, $\phi_m(Z_m)$ converges to 0 for any sequence $Z_m$ of $D[0,1]$ that converges to $Z \in C[0,1]$. Therefore, according to the Extended Continuous Mapping Theorem [27, Theorem 18.11], $\phi_m(\mathbb{Z}_m)$ converges in distribution (hence also in probability) to 0. □

*Proof of Proposition 3.6.* For $m \in \mathbb{N}$, let $\mathcal{T}_m^\varepsilon : D[0,1] \to [0,1]$ be defined by

$$\forall F \in D[0,1], \mathcal{T}_m^\varepsilon(F) = \mathcal{T}\left(F - \varepsilon_m\right)\,.$$

Let $\tau^\star = \mathcal{T}(G)$, $\widehat{\tau} = \mathcal{T}\big(\widehat{\mathbb{G}}_m\big)$, $\widehat{\tau}_m = \mathcal{T}_m\big(\widehat{\mathbb{G}}_m\big)$ and $\widehat{\tau}_m^\varepsilon = \mathcal{T}_m^\varepsilon\big(\widehat{\mathbb{G}}_m\big)$. Write

$$\mathsf{FDP}_m(t) = \frac{\pi_0 \widehat{\mathbb{G}}_{0,m}(t)}{\widehat{\mathbb{G}}_m(t)},$$

where $\widehat{\mathbb{G}}_{0,m}$ and $\widehat{\mathbb{G}}_m$ are non decreasing functions, so that

$$\frac{\pi_0 \widehat{\mathbb{G}}_{0,m}\left(\widehat{\tau}_m^\varepsilon\right)}{\widehat{\mathbb{G}}_m\left(\widehat{\tau}\right)} \le \mathsf{FDP}_m(\widehat{\tau}_m) \le \frac{\pi_0 \widehat{\mathbb{G}}_{0,m}\left(\widehat{\tau}\right)}{\widehat{\mathbb{G}}_m\left(\widehat{\tau}_m^\varepsilon\right)}$$

because $\widehat{\tau}_m^\varepsilon \le \widehat{\tau}_m \le \widehat{\tau}$. Therefore,

$$\begin{aligned}\frac{\pi_0\left(\widehat{\mathbb{G}}_{0,m}\left(\widehat{\tau}_m^\varepsilon\right) - \widehat{\mathbb{G}}_{0,m}\left(\widehat{\tau}\right)\right)}{\widehat{\mathbb{G}}_m\left(\widehat{\tau}\right)} &\le \mathsf{FDP}_m(\widehat{\tau}_m) - \mathsf{FDP}_m(\widehat{\tau}) \\ \le \frac{\pi_0\left(\widehat{\mathbb{G}}_{0,m}\left(\widehat{\tau}\right) - \widehat{\mathbb{G}}_{0,m}\left(\widehat{\tau}_m^\varepsilon\right)\right)}{\widehat{\mathbb{G}}_m\left(\widehat{\tau}_m^\varepsilon\right)} &- \left(\mathsf{FDP}_m(\widehat{\tau}) - \mathsf{FDP}_m(\widehat{\tau}_m^\varepsilon)\right)\end{aligned}$$

As $\mathcal{T}_m^\varepsilon \leq \mathcal{T}_m \leq \mathcal{T}$ for any $m \in \mathbb{N}$, Lemma 7.4 ensures that $\sqrt{m}\left(\widehat{\tau}_m - \widehat{\tau}\right)$ and $\sqrt{m}\left(\widehat{\tau}_m^\varepsilon - \widehat{\tau}\right)$ converge to 0 in probability. Therefore, as $\widehat{\mathbb{G}}_m\left(\widehat{\tau}\right)$ and $\widehat{\mathbb{G}}_m\left(\widehat{\tau}_m^\varepsilon\right)$ converge in probability to $G\left(\mathcal{T}(G)\right)$ as $m \to +\infty$, we have

$$\frac{\pi_0\sqrt{m}\left(\widehat{\mathbb{G}}_{0,m}\left(\widehat{\tau}_m^\varepsilon\right) - \widehat{\mathbb{G}}_{0,m}\left(\widehat{\tau}\right)\right)}{\widehat{\mathbb{G}}_m\left(\widehat{\tau}\right)} \xrightarrow{P} 0$$

and

$$\frac{\pi_0\sqrt{m}\left(\widehat{\mathbb{G}}_{0,m}\left(\widehat{\tau}\right) - \widehat{\mathbb{G}}_{0,m}\left(\widehat{\tau}_m^\varepsilon\right)\right)}{\widehat{\mathbb{G}}_m\left(\widehat{\tau}_m^\varepsilon\right)} \xrightarrow{P} 0,$$

which concludes the proof because $\sqrt{m}\left(\mathsf{FDP}_m(\widehat{\tau}) - \mathsf{FDP}_m(\widehat{\tau}_m^\varepsilon)\right)$ also converges in probability to 0 (according to Lemma 7.4). □

### 7.2. *Asymptotic* FDP*: specific threshold functions*

We now apply the results of section 7.1 to threshold functions of the form

$$\mathcal{T}(F) = \mathcal{U}\left(F, \mathcal{A}(F)\right),$$

with

$$\mathcal{U}(F,\alpha) = \sup\{u \in [0,1], F(u) \geq r_\alpha(u)\}.$$

In this section, we will use the notation $r : (\alpha, u) \mapsto r_\alpha(u)$ whenever the dependence of $r_\alpha$ in $\alpha$ is of importance. We begin by giving sufficient conditions for the regularity of $\mathcal{U}$ and $\mathcal{A}$ under which Condition C.1 holds (section 7.2.1, Proposition 7.5). Then we provide sufficient conditions for $\mathcal{U}$ to be regular enough to be consistent with hypotheses $(i)$ to $(iii)$ of Proposition 7.5 (section 7.2.2). Finally we derive the asymptotic distribution of the corresponding False Discovery Proportion (section 7.2.3).

#### *7.2.1. Hadamard differentiability of $\mathcal{T}$*

**Proposition 7.5** (Hadamard differentiability of $\mathcal{T}$)**.** *Let $C[0,1]$ be the set of continuous functions of $D[0,1]$. Suppose that*

- *(i) $\mathcal{U}$ is Hadamard-differentiable with respect to its first variable at $(G,\alpha)$, tangentially to $C[0,1]$, for any $\alpha$ in a neighborhood of $\mathcal{A}(G)$; its derivative will be denoted by $\nabla_F\mathcal{U}_{G,\alpha}$;*
- *(ii) $\nabla_F\mathcal{U}_{G,\cdot}$ is continuous at $\mathcal{A}(G)$;*
- *(iii) $\mathcal{U}$ is differentiable with respect to its second variable; its derivative will be denoted by $\nabla_\alpha\mathcal{U}_{(G,\mathcal{A}(G))}$;*
- *(iv) $\mathcal{A}$ is Hadamard-differentiable at $G$ tangentially to $C[0,1]$; its derivative will be denoted by $\dot{\mathcal{A}}_G$.*

*Then, $\mathcal{T}$ is Hadamard-differentiable at $G$ tangentially to $C[0,1]$, with derivative $\dot{\mathcal{T}}_G$, defined for any $H \in D[0,1]$ by*

$$\dot{\mathcal{T}}_G(H) = \nabla_F \mathcal{U}_{(G,\mathcal{A}(G))}(H) + \dot{\mathcal{A}}_G(H) \nabla_\alpha \mathcal{U}_{(G,\mathcal{A}(G))} .$$

*Proof of Proposition 7.5.* Let $H \in C[0,1]$, and $H_t$ be a family of functions of $D[0,1]$ that converges to $H$ on $(D[0,1], \|.\|_\infty)$ as $t \to 0$.

As $\mathcal{A}$ is continuous at $G$ (by $(iv)$), $\mathcal{A}(G + tH_t)$ lies in a neighborhood of $\mathcal{A}(G)$ for small $t > 0$, and $\mathcal{U}$ is Hadamard-differentiable with respect to its first variable at $(G, \mathcal{A}(G + tH_t))$ by $(i)$. We therefore have

$$\begin{aligned}
\mathcal{T}(G + tH_t) &= \mathcal{U}(G + tH_t, \mathcal{A}(G + tH_t)) \\
&= \mathcal{U}(G, \mathcal{A}(G + tH_t)) + t \nabla_F \mathcal{U}_{G,\mathcal{A}(G+tH_t)}(H)(1 + o(1)) \\
&= \mathcal{U}(G, \mathcal{A}(G + tH_t)) + t \nabla_F \mathcal{U}_{G,\mathcal{A}(G)}(H)(1 + o(1))
\end{aligned}$$

by the continuity of $\nabla_F \mathcal{U}_{G,\cdot}$ at $\mathcal{A}(G)$ $(ii)$. Then, combining $(iii)$ and $(iv)$ yields

$$\begin{aligned}
\mathcal{U}(G, \mathcal{A}(G + tH_t)) &= \mathcal{U}\left(G, \mathcal{A}(G) + t\dot{\mathcal{A}}_G(H)(1 + o(1))\right) \\
&= \mathcal{U}(G, \mathcal{A}(G)) + t\dot{\mathcal{A}}_G(H) \nabla_\alpha \mathcal{U}_{(G,\mathcal{A}(G))} (1 + o(1)) \\
&= \mathcal{T}(G) + t\dot{\mathcal{A}}_G(H) \nabla_\alpha \mathcal{U}_{(G,\mathcal{A}(G))} (1 + o(1))
\end{aligned}$$

so that

$$\lim_{t \to 0} \frac{\mathcal{T}(G + tH_t) - \mathcal{T}(G)}{t} = \nabla_F \mathcal{U}_{(G,\mathcal{A}(G))}(H) \dot{\mathcal{A}}_G(H) \nabla_\alpha \mathcal{U}_{(G,\mathcal{A}(G))},$$

which concludes the proof. □

#### 7.2.2. *Regularity of $\mathcal{U}$*

The crucial point for proving the desired regularity of $\mathcal{U}$ is its Hadamard differentiability with respect to its first variable at $(G, \alpha)$, tangentially to $C[0,1]$, for $\alpha$ in a neighborhood of $\mathcal{A}(G)$. Lemma 7.6 is a straightforward analytical translation of Conditions C.2 and C.3.

**Lemma 7.6.** *Under Conditions C.2 and C.3, the unique interior right crossing point $\tau^\star$ between $r_\alpha$ and $G$ is positive. If $r$ is $C^1$ on $(0,1] \times [0,1]$, there exists a neighborhood $V = A \times U$ of $(\mathcal{A}(G), \tau^\star)$ such that for any $(\alpha, x) \in V$, $\psi_{G,\alpha}$ : $u \mapsto r_\alpha(u) - G(u)$ is locally inversible around $\mathcal{U}(G, \alpha)$, with $\dot{\psi}_{G,\alpha}(x) > 0$.*

We begin by proving the continuity of $\mathcal{U}$ at $(G, \alpha)$ for $\alpha$ in a neighborhood of $\mathcal{A}(G)$. We then (Proposition 7.11) provide the sufficient conditions for conditions $(i)$, $(ii)$, and $(iii)$ of Proposition 7.5 to hold.

**Lemma 7.7.** *For any $F \in D[0,1]$, and $\alpha \in [0,1]$ such that $r_\alpha$ is continuous, one of the following two assertions holds:*

*(i)* $F(\mathcal{U}(F, \alpha)) = r_\alpha(\mathcal{U}(F, \alpha))$

*(ii)* $F(\mathcal{U}(F,\alpha)) \le r_\alpha(\mathcal{U}(F,\alpha)) \le F(\mathcal{U}(F,\alpha)^-)$

*Proof of Lemma 7.7.* According to the definition of $\mathcal{U}(F,\alpha)$, $F(u) \le r_\alpha(u)$ for any $u > \mathcal{U}(F,\alpha)$. As $F$ is right-continuous and $r_\alpha$ is continuous, we have $F(\mathcal{U}(F,\alpha)) \le r_\alpha(\mathcal{U}(F,\alpha))$. Therefore, either $(i)$ holds, or $F(\mathcal{U}(F,\alpha)) < r_\alpha(\mathcal{U}(F,\alpha))$. In the second case, according to the definition of $\mathcal{U}(F,\alpha)$, there is a non decreasing sequence $(u_n)$ that converges to $\mathcal{U}(F,\alpha)$ such that $F(u_n) \ge r_\alpha(u_n)$. As $r_\alpha$ is continuous and $F$ is left-continuous, we have $F(\mathcal{U}(F,\alpha)^-) \ge r_\alpha(\mathcal{U}(F,\alpha))$, which proves $(ii)$. □

**Proposition 7.8.** *Let $F \in D[0,1]$ be non decreasing, and $\alpha = \mathcal{A}(F)$. If $r_\alpha$ is continuous, then*

$$F(\mathcal{U}(F,\alpha)) = r_\alpha(\mathcal{U}(F,\alpha)) .$$

*Proof of Proposition 7.8.* Let us consider the two assertions of Lemma 7.7: as $F$ is non decreasing, $(ii)$ can be reduced to $(i)$. □

**Proposition 7.9.** *Let $r$ be continuous on $(0,1] \times [0,1]$. Let $F \in D[0,1]$ be non decreasing, and $\alpha \in (0,1]$. Let $F_t$ be a sequence of functions of $D[0,1]$ such that $(F_t)_{t>0}$ converges to $F$ on $(D[0,1], \|.\|_\infty)$ as $t \to 0$, and $\alpha_t \to \alpha$ as $t \to 0$. Denote by $\psi_{F,\alpha}$ the function defined on $[0,1]$ by*

$$\forall u \in [0,1], \psi_{F,\alpha}(u) = r_\alpha(u) - F(u) .$$

*Then,*

$$\lim_{t\to 0} \psi_{F,\alpha}\left(\mathcal{U}(F_t,\alpha_t)\right) = \psi_{F,\alpha}\left(\mathcal{U}(F,\alpha)\right) .$$

*Proof of Proposition 7.9.* For each fixed $t \in [0,1]$, one of the following two assertions holds according to Lemma 7.7:

(i) $F_t(\mathcal{U}(F_t,\alpha_t)) = r_{\alpha_t}(\mathcal{U}(F_t,\alpha_t))$
(ii) $F_t(\mathcal{U}(F_t,\alpha_t)) \le r_{\alpha_t}(\mathcal{U}(F_t,\alpha_t)) \le F_t(\mathcal{U}(F_t,\alpha_t)^-)$.

If $(ii)$ holds, then, as $F$ is non decreasing we have

$$(F - F_t)(\mathcal{U}(F_t,\alpha_t)^-) \le F(\mathcal{U}(F_t,\alpha_t)) - r_{\alpha_t}(\mathcal{U}(F_t,\alpha_t)) \le (F - F_t)(\mathcal{U}(F_t,\alpha_t)) .$$

If $(i)$ holds, then $F(\mathcal{U}(F_t,\alpha_t)) - r_{\alpha_t}(\mathcal{U}(F_t,\alpha_t)) = F(\mathcal{U}(F_t,\alpha_t)) - F_t(\mathcal{U}(F_t,\alpha_t))$. In either case, we have $|F(\mathcal{U}(F_t,\alpha_t)) - r_{\alpha_t}(\mathcal{U}(F_t,\alpha_t))| \le \|F - F_t\|$, which tends to 0 as $t \to 0$. As $r$ is continuous on $(0,1] \times [0,1]$, $r_{\alpha_t}$ converges uniformly to $r_\alpha$ on the compact $[\alpha/2, 1]$, and we have

$$\lim_{t\to 0} F(\mathcal{U}(F_t,\alpha_t)) - r_\alpha(\mathcal{U}(F_t,\alpha_t)) = 0 ,$$

which concludes the proof as $\psi_{F,\alpha}\left(\mathcal{U}(F,\alpha)\right) = 0$. □

**Corollary 7.10** (Continuity of $\mathcal{U}$)**.** *Let $r$ be $C^1$ on $(0,1] \times [0,1]$. According to Conditions C.2 and C.3, there is a neighborhood $A$ of $\mathcal{A}(G)$ such that $\mathcal{U}$ is continuous at $(G,\alpha)$ for any $\alpha \in A$.*

**Proposition 7.11** (Differentiability of $\mathcal{U}$)**.** *Let us assume that $r$ is $C^1$ on $(0,1]\times[0,1]$. Under Conditions C.2 and C.3,*

*(i) $\mathcal{U}$ is Hadamard-differentiable with respect to its first variable at $(G,\alpha)$, tangentially to $C[0,1]$ for any $\alpha$ in a neighborhood $A$ of $\mathcal{A}(G)$, with derivative $\nabla_F\mathcal{U}_{G,\alpha}$ defined for any $H\in C[0,1]$ by*

$$\nabla_F\mathcal{U}_{G,\alpha}(H)=\frac{H(\mathcal{U}(G,\alpha))}{\frac{\partial r}{\partial u}(\alpha,\mathcal{U}(G,\alpha))-g(\mathcal{U}(G,\alpha))}$$

*(ii) $\nabla_F\mathcal{U}_{G,.}$ is continuous at $\mathcal{A}(G)$ on $A$.*

*(iii) $\mathcal{U}$ is differentiable with respect to its second variable, with derivative*

$$\nabla_\alpha\mathcal{U}_{G,\mathcal{A}(G)}=-\frac{\frac{\partial r}{\partial \alpha}(\mathcal{A}(G),\tau^\star))}{\frac{\partial r}{\partial u}(\mathcal{A}(G),\tau^\star))-g(\tau^\star))},$$

*where $\tau^\star=\mathcal{U}(G,\mathcal{A}(G))$.*

**Corollary 7.12.** *If we also assume that $\mathcal{A}$ is Hadamard-differentiable at $G$, tangentially to $C[0,1]$, then $\mathcal{T}$ is Hadamard-differentiable at $G$, tangentially to $C[0,1]$, with derivative defined for any $H\in C[0,1]$ by*

$$\dot{\mathcal{T}}_G(H)=\frac{H(\tau^\star)-\frac{\partial r}{\partial \alpha}(\mathcal{A}(G),\tau^\star))\dot{\mathcal{A}}_G(H)}{\frac{\partial r}{\partial u}(\mathcal{A}(G),t)-g(t)}$$

*Proof of Proposition 7.11.* Let $\tau^\star=\mathcal{U}(G,\mathcal{A}(G))$. Throughout the proof, $V=A\times U$ denotes the neighborhood of $(\mathcal{A}(G),\tau^\star)$ defined in Lemma 7.6.

(i) Let $\alpha\in A$. Let $H\in C[0,1]$, and $H_t$ be a family of functions of $D[0,1]$ that converges to $H$ on $(D[0,1],\|.\|_\infty)$ as $t\to 0$. Let $v=\mathcal{U}(G,\alpha)$ and $v_t=\mathcal{U}(G_t,\alpha)$, with $G_t=G+tH_t$. By the continuity of $U$ (Corollary 7.10), $v_t\to v$ as $t\to 0$. Therefore, applying Taylor's formula to $\psi_{G,\alpha}:u\mapsto r_\alpha(u)-G(u)$ yields

$$\psi_{G,\alpha}(v_t)-\psi_{G,\alpha}(v)\underset{t\to 0}{=}(v_t-v)\dot{\psi}_{G,\alpha}(v)\left(1+o(1)\right).$$

As $\alpha\in A$, we have $\dot{\psi}_{G,\alpha}(v)=\frac{\partial r}{\partial u}(\alpha,v)-g(v)>0$. Therefore, as $\psi_{G,\alpha}(v)=0$ according to Proposition 7.8,

$$v_t-v\underset{t\to 0}{=}\frac{\psi_{G,\alpha}(v_t)}{\frac{\partial r}{\partial u}(\alpha,v)-g(v)}\left(1+o(1)\right)$$

so that it is sufficient to prove that $\lim_{t\to 0}\psi_{G,\alpha}(v_t)/t=H(v)$. The behavior of $\psi_{G,\alpha}(v_t)=r_\alpha(v_t)-G(v_t)$ can be derived using the same argument as in the proof of proposition 7.9; Lemma 7.7, we have either $G_t(v_t)=r_\alpha(v_t)$ or $G_t(v_t)\le r_\alpha(v_t)\le G_t(v_t^-)$. In the first case, $r_\alpha(v_t)-G(v_t)=(G_t-G)(v_t)=tH_t(v_t)$, and $\lim_{t\to 0}\frac{r_\alpha(v_t)-G(v_t)}{t}=H(t)$ according to Lemma 7.1. In the second case, we have

$$(G_t-G)(v_t)\le r_\alpha(v_t)-G(v_t)\le(G_t-G)(v_t^-)$$

as $G$ is non decreasing, that is, $tH_t(v_t) \leq \psi_{G,\alpha}(v_t) \leq tH_t(v_t^-)$. Therefore, we have

$$H_t(v_t) - H(v_t) \leq \frac{\psi_{G,\alpha}(v_t)}{t} - H(v_t) \leq H_t(v_t^-) - H(v_t).$$

As $H$ is continuous, $H(v_t) = H(v_t^-)$, and the upper and lower bounds converge to 0 according to Lemma 7.1, and $(i)$ is proved.

(ii) is a consequence of the continuity of $\mathcal{U}$ (with respect to its second variable), that of $g$ and that of $\nabla_u r$ with respect to its first variable.

(iii) Let $\alpha = \mathcal{A}(G)$. Let $\alpha_t$ be a sequence of points of $(0, 1]$ that converges to $\alpha$. Let $v = \mathcal{U}(G, \alpha)$ and $v_t = \mathcal{U}(G, \alpha_t)$. By the continuity of $U$ (Proposition 7.9), $v_t \to v$ as $t \to 0$. Therefore, applying Taylor's formula to $\psi_{G,\alpha} : u \mapsto r_\alpha(u) - G(u)$ yields

$$\psi_{G,\alpha}(v_t) - \psi_{G,\alpha}(v) \underset{t\to 0}{=} (v_t - v)\dot{\psi}_{G,\alpha}(v)\,(1 + o(1)) .$$

We have $\psi_{G,\alpha}(v) = 0$ and $\psi_{G,\alpha}(v_t) = r(\alpha, v_t) - r(\alpha_t, v_t)$ by Proposition 7.8. As $r$ is $C^1$ in a neighborhood of $(\alpha, v)$, we have, according to Taylor's formula,

$$r(\alpha_t, v_t) - r(\alpha, v_t) \underset{t\to 0}{=} (\alpha_t - \alpha)\frac{\partial r}{\partial \alpha}(\alpha, v)\,(1 + o(1)) .$$

As $\alpha = \mathcal{A}$ and $v = \mathcal{U}(G, \alpha)$, $(\alpha, v) \in V$. Therefore $\dot{\psi}_{G,\alpha}(v) = \frac{\partial r}{\partial u}(\alpha, v) - g(v) > 0$ according to Lemma 7.6. Finally, we have

$$\lim_{t\to 0} \frac{r(\alpha, v_t) - r(\alpha_t, v_t)}{\alpha_t - \alpha} = -\frac{\frac{\partial r}{\partial \alpha}(\alpha, v)}{\frac{\partial r}{\partial u}(\alpha, v) - g(v)},$$

which concludes the proof. □

### 7.2.3. *Asymptotic* FDP

**Theorem 7.13** (Asymptotic distribution of $\mathsf{FDP}_m$). *Let $r_\alpha$ be a rejection curve such that $r$ is $C^1$ on $(0, 1] \times [0, 1]$, and $\mathcal{A}$ a level function. Let us denote by $\mathcal{T} : F \mapsto \mathcal{U}(F, \mathcal{A}(F))$ the associated threshold function, where $\mathcal{U}(F, \alpha) = \sup\{u \in [0, 1], F(u) \geq r_\alpha(u)\}$.*

*Under Conditions C.2 and C.3, if $\mathcal{A}$ is Hadamard-differentiable at $G$ tangentially to $C[0, 1]$, then*

*(i)*

$$\sqrt{m}\left(\mathcal{T}(\widehat{\mathbb{G}}_m) - \tau^\star\right) \rightsquigarrow \frac{\mathbb{Z}(\tau^\star) - \frac{\partial r}{\partial \alpha}(\mathcal{A}(G), \tau^\star))\dot{\mathcal{A}}_G(\mathbb{Z})}{\frac{\partial r}{\partial u}(\mathcal{A}(G), t) - g(t)},$$

*(ii)*

$$\sqrt{m}\left(\mathsf{FDP}_m(\mathcal{T}(\widehat{\mathbb{G}}_m)) - \mathsf{p}^\star\right) \rightsquigarrow X,$$

*where* $\mathsf{p}^\star = \frac{\pi_0\tau^\star}{G(\tau^\star)}$ *is the* pFDR *achieved by procedure* $\mathcal{T}$, $\dot{\mathsf{p}}(t) = \frac{\pi_0}{G(t)}\left(1 - \frac{tg(t)}{G(t)}\right)$, *and*

$$X = \mathsf{p}^\star(1-\mathsf{p}^\star\zeta(\tau^\star))\frac{\mathbb{Z}_0(\tau^\star)}{\tau^\star} + \mathsf{p}^\star(1-\mathsf{p}^\star)\zeta(\tau^\star)\frac{\mathbb{Z}_1(\tau^\star)}{G_1(\tau^\star)} + \dot{\mathsf{p}}(\tau^\star)\xi(\tau^\star)\dot{\mathcal{A}}_G(\mathbb{Z})\,,$$

*with* $\zeta(t) = -\frac{\frac{G(t)}{t} - \frac{\partial r}{\partial u}(\mathcal{A}(G),t)}{\frac{\partial r}{\partial u}(\mathcal{A}(G),t)-g(t)}$, $\xi(t) = \frac{-\frac{\partial r}{\partial \alpha}(\mathcal{A}(G),t)}{\frac{\partial r}{\partial u}(\mathcal{A}(G),t)-g(t)}$, *and* $\mathbb{Z} = \pi_0\mathbb{Z}_0 + (1-\pi_0)\mathbb{Z}_1$, *where* $\mathbb{Z}_0$ *and* $\mathbb{Z}_1$ *are independent Gaussian processes such that* $\mathbb{Z}_0 \overset{(d)}{=} \mathbb{B}$ *and* $\mathbb{Z}_1 \overset{(d)}{=} \mathbb{B}\circ G_1$, *and* $\mathbb{B}$ *is a standard Brownian bridge on* $[0,1]$.

*Proof of Theorem 7.13.* Under these assumptions, Condition C.1 holds for $\mathcal{T}$ according to Corollary 7.12, with

$$\dot{\mathcal{T}}_G(H) = \frac{H(\tau^\star) - \frac{\partial r}{\partial \alpha}(\mathcal{A}(G),\tau^\star))\dot{\mathcal{A}}_G(H)}{\frac{\partial r}{\partial u}(\mathcal{A}(G),t) - g(t)}$$

Therefore, Theorem 3.2 yields $\sqrt{m}\left(\mathcal{T}(\widehat{\mathbb{G}}_m) - \tau^\star)\right) \rightsquigarrow \dot{\mathcal{T}}_G(\mathbb{Z})$, and

$$\sqrt{m}\left(\mathsf{FDP}_m(\mathcal{T}(\widehat{\mathbb{G}}_m)) - \frac{\pi_0\tau^\star}{G(\tau^\star)}\right) \rightsquigarrow X\,,$$

with $X = \mathsf{p}^\star(1-\mathsf{p}^\star)\left(\frac{\mathbb{Z}_0(\tau^\star)}{\tau^\star} - \frac{\mathbb{Z}_1(\tau^\star)}{G_1(\tau^\star)}\right) + \dot{\mathsf{p}}(\tau^\star)\dot{\mathcal{T}}_G(\mathbb{Z})$ and $\mathbb{Z} = \pi_0\mathbb{Z}_0 + (1-\pi_0)\mathbb{Z}_1$, where $\mathbb{Z}_0$ and $\mathbb{Z}_1$ are independent Gaussian processes such that $\mathbb{Z}_0 \overset{(d)}{=} \mathbb{B}$ and $\mathbb{Z}_1 \overset{(d)}{=} \mathbb{B}\circ G_1$, and $\mathbb{B}$ is a standard Brownian bridge on $[0,1]$.

Letting $\delta(t) = \frac{1}{\frac{\partial r}{\partial u}(\mathcal{A}(G),t)-g(t)}$, we have

$$\dot{\mathcal{T}}_G(\mathbb{Z}) = \delta(\tau^\star)\left(\mathbb{Z}(\tau^\star) - \frac{\partial r}{\partial \alpha}(\mathcal{A}(G),\tau^\star))\dot{\mathcal{A}}_G(H)\right)\,.$$

As $\dot{\mathsf{p}}(t) = \mathsf{p}(t)\left(\frac{1}{t} - \frac{g(t)}{G(t)}\right)$, we have

$$\dot{\mathsf{p}}(\tau^\star)\delta(\tau^\star)\mathbb{Z}(\tau^\star) = \mathsf{p}^\star\left(\frac{G(\tau^\star)}{\tau^\star} - g(\tau^\star)\right)\delta(\tau^\star)\frac{\mathbb{Z}(\tau^\star)}{G(\tau^\star)}\,,$$

with $\frac{\mathbb{Z}(\tau^\star)}{G(\tau^\star)} = \mathsf{p}^\star\frac{\mathbb{Z}_0(\tau^\star)}{\tau^\star} + (1-\mathsf{p}^\star)\frac{\mathbb{Z}_1(\tau^\star)}{G_1(\tau^\star)}$. Hence letting $\zeta(t) = 1-\delta(t)\left(\frac{G(t)}{t} - g(t)\right)$, we have

$$X = \mathsf{p}^\star(1-\mathsf{p}^\star\zeta(\tau^\star))\frac{\mathbb{Z}_0(\tau^\star)}{\tau^\star} - \mathsf{p}^\star(1-\mathsf{p}^\star)\zeta(\tau^\star)\frac{\mathbb{Z}_1(\tau^\star)}{G_1(\tau^\star)} + \dot{\mathsf{p}}(\tau^\star)\xi(\tau^\star)\dot{\mathcal{A}}_G(\mathbb{Z})\,,$$

where $\xi(t) = \frac{-\frac{\partial r}{\partial \alpha}(\mathcal{A}(G),t)}{\frac{\partial r}{\partial u}(\mathcal{A}(G),t)-g(t)}$. This concludes the proof since $\zeta$ may be written as $\zeta(t) = -\frac{\frac{G(t)}{t} - \frac{\partial r}{\partial u}(\mathcal{A}(G),t)}{\frac{\partial r}{\partial u}(\mathcal{A}(G),t)-g(t)}$. □

### 7.3. Limit distribution for procedures under consideration

#### 7.3.1. One-stage procedures

In this section $\mathcal{A}(G)$ is fixed. Therefore, only the dependence of $r_\alpha$ $u$ is of importance. In order to lighten the notation we let

$$\dot{r_\alpha} = \frac{\partial r}{\partial u}(\alpha, \cdot)\,.$$

**Theorem 7.14** (Asymptotic FDP for one-stage procedures)**.** *Let $\mathcal{T} : F \mapsto \mathcal{U}(F,\alpha)$ a one-stage procedure such that $r_\alpha$ is continuous on $[0,1]$, and $C^1$ in a neighborhood of $\tau^\star = \mathcal{T}(G)$. Under Condition C.2 and C.3,*

*(i)*

$$\sqrt{m}\left(\mathcal{T}(\widehat{\mathbb{G}}_m) - \tau^\star)\right) \rightsquigarrow \frac{\mathbb{Z}(\tau^\star)}{\dot{r_\alpha}(\tau^\star) - g(\tau^\star)}$$

*(ii)*

$$\sqrt{m}\left(\mathsf{FDP}_m(\mathcal{T}(\widehat{\mathbb{G}}_m)) - \mathsf{p}^\star\right) \rightsquigarrow X\,,$$

*with*

$$X = \mathsf{p}^\star(1-\mathsf{p}^\star\zeta(\tau^\star))\frac{\mathbb{Z}_0(\tau^\star)}{\tau^\star} - \mathsf{p}^\star(1-\mathsf{p}^\star))\zeta(\tau^\star)\frac{\mathbb{Z}_1(\tau^\star)}{G_1(\tau^\star)}\,,$$

*where $\mathsf{p}^\star = \frac{\pi_0\tau^\star}{r_\alpha(\tau^\star)}$, $\zeta(t) = \frac{\frac{\partial r}{\partial u}(\mathcal{A}(G),t) - \frac{G(t)}{t}}{\frac{\partial r}{\partial u}(\mathcal{A}(G),t) - g(t)}$, and $\mathbb{Z} = \pi_0\mathbb{Z}_0 + (1-\pi_0)\mathbb{Z}_1$, where $\mathbb{Z}_0$ and $\mathbb{Z}_1$ are independent Gaussian processes such that $\mathbb{Z}_0 \overset{(d)}{=} \mathbb{B}$ and $\mathbb{Z}_1 \overset{(d)}{=} \mathbb{B}\circ G_1$, and $\mathbb{B}$ is a standard Brownian bridge on $[0,1]$.*

*Proof of Theorem 7.14.* As $\mathcal{T}$ is a one-stage procedure, we have $\mathcal{A} = \alpha$. Therefore, the assumptions for Theorem 7.13 hold, with $\xi = 0$. According to Proposition 7.8, $G(\tau^\star) = r_\alpha(\tau^\star)$; therefore $\mathsf{p}^\star = \frac{\pi_0\tau^\star}{r_\alpha(\tau^\star)}$, which concludes the proof. □

*Proof of Theorem 4.2* (BH95). Uniqueness Condition C.3 always holds because $r_\alpha$ is linear, and Condition C.2 holds because it corresponds to Condition C.4. Therefore, Theorem 7.14 can be applied, and we have $\zeta(\tau^\star) = 0$ since $\dot{r_\alpha}(\tau^\star) = 1/\alpha = r_\alpha(\tau^\star)/\tau^\star$, and $\mathsf{p}(\tau^\star) = \pi_0\alpha$. Hence,

$$\sqrt{m}\left(\mathcal{T}(\widehat{\mathbb{G}}_m) - \tau^\star)\right) \rightsquigarrow \frac{\mathbb{Z}(\tau^\star)}{1/\alpha - g(\tau^\star)}$$

and

$$\sqrt{m}\left(\mathsf{FDP}_m(\mathcal{T}(\widehat{\mathbb{G}}_m)) - \pi_0\alpha\right) \rightsquigarrow \pi_0\alpha\frac{\mathbb{Z}_0(\tau^\star)}{\tau^\star},$$

which concludes the proof because $\operatorname{Var}\mathbb{Z}_0(\tau^\star) = \tau^\star(1-\tau^\star)$. □

*Proof of Theorem 4.7* (FDR08). The uniqueness Condition C.3, and existence Conditions C.4 and C.7 ensure that there is a unique interior right crossing point $\tau^\star$ between $f_\alpha$ and $G$, which satisfies $\tau^\star \leq \kappa$. Condition C.6 guarantees that $\tau^\star$ is also the only right crossing point between $f_\alpha^\lambda$ and $G$. Thus, $[0, \kappa]$ is a neighborhood of $\tau^\star$ in which $f_\alpha^\lambda$ coincides with $f_\alpha$ and is $C^1$, with $\dot{f}_\alpha(u) = \frac{\alpha}{(\alpha+(1-\alpha)u)^2}$. Therefore, Theorem 7.14 yields $\sqrt{m}\big(\mathsf{FDP}_m(\mathcal{T}(\widehat{\mathbb{G}}_m)) - \mathsf{p}^\star\big) \rightsquigarrow X$, with

$$X = \mathsf{p}^\star(1 - \mathsf{p}^\star\zeta(\tau^\star))\frac{\mathbb{Z}_0(\tau^\star)}{\tau^\star} - \mathsf{p}^\star(1 - \mathsf{p}^\star)\zeta(\tau^\star)\frac{\mathbb{Z}_1(\tau^\star)}{G_1(\tau^\star)},$$

where $\mathsf{p}^\star = \frac{\pi_0\tau^\star}{f_\alpha(\tau^\star)} = \pi_0(\alpha + (1-\alpha)\tau^\star)$ and $\zeta(\tau^\star) = -\frac{G(\tau^\star)/\tau^\star - \dot{f}_\alpha(\tau^\star)}{\dot{f}_\alpha(\tau^\star) - g(\tau^\star)}$. Letting

$$\overline{\pi_0}(t) = \frac{1 - G(t)}{1 - t},$$

we have $G(\tau^\star)/\tau^\star = \overline{\pi_0}(\tau^\star)/\alpha$, and $\dot{f}_\alpha(\tau^\star) = \alpha(f_\alpha(\tau^\star)/\tau^\star)^2 = \overline{\pi_0}(\tau^\star)^2/\alpha$, so that $\mathsf{p}^\star = \alpha\pi_0/\overline{\pi_0}(\tau^\star)$, and

$$\begin{aligned}
\zeta(\tau^\star) &= -\frac{G(\tau^\star)/\tau^\star - \dot{f}_\alpha(\tau^\star)}{\dot{f}_\alpha(\tau^\star) - g(\tau^\star)} \\
&= -\frac{\overline{\pi_0}(\tau^\star)/\alpha - \overline{\pi_0}(\tau^\star)^2/\alpha}{\overline{\pi_0}(\tau^\star)^2/\alpha - g(\tau^\star)} \\
&= -(1 - \overline{\pi_0}(\tau^\star))\frac{\overline{\pi_0}(\tau^\star)/\alpha}{\overline{\pi_0}(\tau^\star)^2/\alpha - g(\tau^\star)},
\end{aligned}$$

which concludes the proof. □

*Proof of Theorem 4.10* (BR08($\lambda$)). The uniqueness Condition C.3, and existence Conditions C.8 and C.9 ensure that there is a unique interior right crossing point $\tau^\star$ between $b_\alpha^\lambda$ and $G$, which satisfies $\tau^\star \leq \lambda$. Thus $[0, \lambda]$ is a neighborhood of $\tau^\star$ in which $b_\alpha^\lambda$ is $C^1$, with $\dot{b}_\alpha^\lambda(u) = \frac{\alpha(1-\lambda)}{(\alpha(1-\lambda)+u)^2}$. Therefore, Theorem 7.14 yields $\sqrt{m}\big(\mathsf{FDP}_m(\mathcal{T}(\widehat{\mathbb{G}}_m)) - \mathsf{p}^\star\big) \rightsquigarrow X$, with

$$X = \mathsf{p}^\star(1 - \mathsf{p}^\star\zeta(\tau^\star))\frac{\mathbb{Z}_0(\tau^\star)}{\tau^\star} - \mathsf{p}^\star(1 - \mathsf{p}^\star)\zeta(\tau^\star)\frac{\mathbb{Z}_1(\tau^\star)}{G_1(\tau^\star)},$$

where $\mathsf{p}^\star = \frac{\pi_0\tau^\star}{b_\alpha(\tau^\star)} = \pi_0(\alpha + (1-\alpha)\tau^\star)$ and $\zeta(\tau^\star) = -\frac{G(\tau^\star)/\tau^\star - \dot{b}_\alpha(\tau^\star)}{\dot{b}_\alpha(\tau^\star) - g(\tau^\star)}$. We have $\dot{f}_\alpha(\tau^\star) = \alpha(1-\lambda)(b_\alpha(\tau^\star)/\tau^\star)^2 = G(\tau^\star)(1 - G(\tau^\star))/\tau^\star$, so that

$$\begin{aligned}
\zeta(\tau^\star) &= -\frac{G(\tau^\star)/\tau^\star - \dot{b}_\alpha(\tau^\star)}{\dot{b}_\alpha(\tau^\star) - g(\tau^\star)} \\
&= -\frac{G(\tau^\star)/\tau^\star(1 - (1 - G(\tau^\star)))}{G(\tau^\star)(1 - G(\tau^\star))/\tau^\star - g(\tau^\star)} \\
&= -\frac{G(\tau^\star)^2/\tau^\star}{G(\tau^\star)(1 - G(\tau^\star))/\tau^\star - g(\tau^\star)},
\end{aligned}$$

which concludes the proof. □

#### 7.3.2. Two-stage adaptive procedures

*Proof of Theorem 4.12.* As pointed out in section 4.3, Condition C.3 always holds because $r_\alpha$ is linear, and Condition C.2 holds as soon as $\mathcal{A}(G) > \alpha^\star$. Therefore, Theorem 7.13 yields $\sqrt{m}\big(\mathsf{FDP}_m(\mathcal{T}(\widehat{\mathbb{G}}_m)) - \mathsf{p}^\star\big) \rightsquigarrow X$, with $\mathsf{p}^\star = \pi_0 \mathcal{A}(G)$, and

$$X = \mathsf{p}^\star(1-\zeta(\tau^\star)\mathsf{p}^\star)\frac{\mathbb{Z}_0(\tau^\star)}{\tau^\star} - \mathsf{p}^\star(1-\mathsf{p}^\star)\zeta(\tau^\star)\frac{\mathbb{Z}(\tau^\star)}{G(\tau^\star)} + \dot{\mathsf{p}}(\tau^\star)\xi(\tau^\star)\dot{\mathcal{A}}_G(\mathbb{Z})\,,$$

where $\dot{\mathsf{p}}(\tau^\star) = \frac{\mathsf{p}^\star}{G(\tau^\star)}\left(\frac{G(\tau^\star)}{\tau^\star} - g(\tau^\star)\right)$, $\zeta(t) = -\frac{\frac{\partial r}{\partial u}(\mathcal{A}(G),t) - G(t)/t}{\frac{\partial r}{\partial u}(\mathcal{A}(G),t) - g(t)}$, and $\xi(t) = \frac{-\frac{\partial r}{\partial \alpha}(\mathcal{A}(G),t)}{\frac{\partial r}{\partial u}(\mathcal{A}(G),t) - g(t)}$. Simes' line is defined by $r_\alpha : u \mapsto u/\alpha$. Therefore, we have $\frac{\partial r}{\partial u}(\mathcal{A}(G),t) = \frac{1}{\mathcal{A}(G)}$ and $\frac{\partial r}{\partial \alpha}(\mathcal{A}(G),t) = -\frac{t}{\mathcal{A}(G)^2}$, and $G(\tau^\star) = \frac{\tau^\star}{\mathcal{A}(G)}$ according to Proposition 7.8. We have $\zeta(\tau^\star) = 0$, $\xi(\tau^\star) = \frac{\tau^\star/\mathcal{A}(G)^2}{\frac{1}{\mathcal{A}(G)} - g(\tau^\star)}$, and $\dot{\mathsf{p}}(\tau^\star) = \mathsf{p}^\star\frac{\mathcal{A}(\mathcal{G})}{\tau^\star}\left(\frac{1}{\mathcal{A}(G)} - g(\tau^\star)\right)$, which concludes the proof. □

**Sto02 procedure.** The following Proposition establishes the Hadamard differentiability of the level function of procedure Sto02. The proof is immediate.

**Proposition 7.15.** *For $F \in D[0,1]$, let*

$$\mathcal{A}(F) = \alpha\frac{1-\lambda}{1-F(\lambda)},$$

*where $\alpha \in [0,1]$. Under Condition C.11, $\mathcal{A}$ is Hadamard-differentiable at $G$, tangentially to $C[0,1]$, with derivative*

$$\dot{\mathcal{A}}_G(H) = \mathcal{A}(G)\frac{H(\lambda)}{1-G(\lambda)}\,.$$

**Proposition 7.16.** *Let $\lambda \in [0,1)$ such that Conditions C.9 and C.11 hold. Then procedures Sto02$(\lambda)$ and STS04$(\lambda)$ are asymptotically equivalent.*

*Proof of Proposition 7.16.* Let $\lambda \in [0,1)$. According to Condition C.9, we have $\mathcal{T}^{\mathsf{Sto02}(\lambda)}(G) < \lambda$. Therefore, procedure Sto02 is asymptotically equivalent to the same procedure truncated at $\lambda$, that is, the procedure with threshold function defined for $F \in D[0,1]$ by

$$\sup\left\{u \in [0,\lambda], F(u) \geq \frac{u}{\alpha}\frac{1-\lambda}{1-F(\lambda)}\right\}\,.$$

We thus work with this truncated version for the remainder of the proof. By definition, the rejection curve of procedure STS04 is larger than that of procedure Sto02. Therefore, we have $\mathcal{T}_m^{\mathsf{STS04}(\lambda)}(F) \leq \mathcal{T}^{\mathsf{Sto02}(\lambda)}(F)$ for any $F \in D[0,1]$. With the same argument we also have $\mathcal{T}^{\mathsf{Sto02}(\lambda)}\left(F - \frac{1}{m}\right) \leq \mathcal{T}_m^{\mathsf{STS04}(\lambda)}(F)$ for any $F \in D[0,1]$. As we have assumed that Condition C.11 holds, Condition C.1 holds for $\mathcal{T}$ according to Proposition 7.15, and the result follows from Proposition 3.6. □

*Proof of Theorem 4.15.* According to Proposition 7.15, and because $\mathcal{A}(G) > \alpha^\star$, Theorem 4.12 ensures that

$$\sqrt{m}\left(\mathsf{FDP}_m(\mathcal{T}^{\mathsf{Sto02}}(\widehat{\mathbb{G}}_m)) - \pi_0\mathcal{A}(G)\right) \rightsquigarrow \pi_0\mathcal{A}(G)\left(\frac{\mathbb{Z}_0(\tau^\star)}{\tau^\star} + \frac{\dot{\mathcal{A}}_G(\mathbb{Z})}{\mathcal{A}(G)}\right),$$

where

$$\dot{\mathcal{A}}_G(\mathbb{Z}) = \mathcal{A}(G)\frac{\mathbb{Z}(\lambda)}{1-G(\lambda)}$$

Denoting $\overline{\pi_0}\,(\lambda) = \overline{\pi_0}^G(\lambda)$, this may be written as

$$\sqrt{m}\left(\mathsf{FDP}_m(\mathcal{T}(\widehat{\mathbb{G}}_m)) - \frac{\pi_0}{\overline{\pi_0}\,(\lambda)}\alpha\right) \rightsquigarrow \frac{\pi_0}{\overline{\pi_0}\,(\lambda)}\alpha\left(\frac{\mathbb{Z}_0(\tau^\star)}{\tau^\star} + \frac{\mathbb{Z}(\lambda)}{1-G(\lambda)}\right)$$

For the calculation of variance, it suffices to note that $\operatorname{Var}\mathbb{Z}_0(\tau^\star) = \tau^\star(1-\tau^\star)$ and

$$\begin{aligned}\mathbb{E}\left[\mathbb{Z}_0(\tau^\star)\mathbb{Z}(\lambda)\right] &= \pi_0\mathbb{E}\left[\mathbb{Z}_0(\tau^\star)\mathbb{Z}_0(\lambda)\right] \\ &= \pi_0\left(\tau^\star \wedge \lambda - \tau^\star\lambda\right),\end{aligned}$$

which concludes the proof since $\mathbb{Z}_0$ and $\mathbb{Z}$ are centered. □

**Procedure** $\mathsf{BKY06}(\lambda)$**[2].** According to Proposition 7.8, we have $F(\mathcal{U}(F,\beta)) = \mathcal{U}(F,\beta)/\beta$ for any $F \in D[0,1]$, so that the level function of procedure BKY06 may be written as

$$\mathcal{A}(F) = \frac{\alpha(1-\lambda)}{1-\mathcal{U}(F,\lambda)/\lambda}.$$

We now prove the Hadamard differentiability of the level function of procedure $\mathsf{BKY06}(\lambda)$ under Condition C.12.

**Proposition 7.17.** *For $\lambda \in [0,1)$ and $F \in D[0,1]$, let $\mathcal{A}(F) = \frac{\alpha(1-\lambda)}{1-\mathcal{U}(F,\lambda)/\lambda}$. Under Condition C.12, $\mathcal{A}$ is Hadamard-differentiable at $G$, tangentially to $C[0,1]$, with derivative*

$$\dot{\mathcal{A}}_G(H) = \frac{\mathcal{A}(G)^2}{\alpha(1-\lambda)}\frac{H(\mathcal{U}(G,\lambda))}{1/(\alpha(1-\lambda)) - g(\mathcal{U}(G,\lambda))}.$$

*Proof of Proposition 7.17.* As Condition C.12 holds, Condition C.4 holds for the BH95 procedure at level $\lambda$: $\mathcal{U}$ is Hadamard-differentiable with respect to its first variable at $(G,\lambda)$, tangentially to $C[0,1]$, with derivative $\nabla_F\mathcal{U}_{G,\lambda}$ defined for any $H \in C[0,1]$ by

$$\nabla_F\mathcal{U}_{G,\lambda}(H) = \frac{H(\mathcal{U}(G,\lambda))}{\frac{\partial r}{\partial u}(\lambda,\mathcal{U}(G,\lambda)) - g(\mathcal{U}(G,\lambda))}.$$

As the rejection curve of $\mathcal{U}$ is Simes' line, we have $\frac{\partial r}{\partial u}(\lambda,\mathcal{U}(G,\lambda)) = \frac{1}{\lambda}$. As $\mathcal{A}(F) = \alpha\phi(\mathcal{U}(F,\lambda))$, where $\phi : x \mapsto (1-\lambda)/(1-x/\lambda)$ is derivable for $x \neq \lambda$, with $\phi'(x) = \frac{\lambda(1-\lambda)}{1-x/\lambda}$, the result follows from the chain rule. □

*Proof of Theorem 4.20.* As Condition C.12 holds, this is a direct consequence of Proposition 7.17 and Theorem 4.12. □

### 7.4. Connections between one-stage and two-stage adaptive procedures

*Proof of Theorem 5.2.* As we have assumed that existence Condition C.2 and uniqueness Condition C.3 hold for procedure $\mathcal{T}$, $\tau^\star = \mathcal{T}(G)$ is the only point in $(0,1)$ such that $G(\tau^\star) = c_\alpha(G,\tau^\star)$. Similarly, as existence Condition C.2 holds for procedure $\mathcal{T}_t$ for any $t \in (0,1)$, $\tau(t) = \mathcal{T}_t(G)$ is the only point in $(0,1)$ such that $G(\tau(t))/\tau(t) = c_\alpha(G,t)/t$. Therefore, we have

$$\begin{aligned} t \leq \tau^\star \quad &\iff \quad G(t) \geq c_\alpha(G,t) \\ &\iff \quad \frac{G(t)}{t} \geq \frac{c_\alpha(G,t)}{t} \\ &\iff \quad \frac{G(t)}{t} \geq \frac{G(\tau(t))}{\tau(t)} \\ &\iff \quad t \leq \tau(t) \end{aligned}$$

as $u \mapsto G(u)/u$, is non increasing (due to the concavity of $G$). As $u \mapsto c_\alpha(G,u)/u$ is non increasing (Condition C.13), we have

$$\begin{aligned} t \leq \tau^\star \quad &\iff \quad \frac{c_\alpha(G,\tau^\star)}{\tau^\star} \leq \frac{c_\alpha(G,t)}{t} \\ &\iff \quad \frac{G(\tau^\star)}{\tau^\star} \leq \frac{G(\tau(t))}{\tau(t)} \\ &\iff \quad \tau(t) \leq \tau^\star\,, \end{aligned}$$

and $(i)$ is proved. Let $\lambda \in (0,1)$. If $\lambda \leq \tau^\star$, then by $(i)$, the sequence $(t_n)$ is non decreasing, and smaller than $\tau^\star$. It therefore converges to a limit $\ell \in [\lambda, \tau^\star]$, such that $\tau(\ell) = \ell$, that is, $G(\ell) = c_\alpha(G,\ell)$. The uniqueness Condition C.3 ensures that $\ell = \tau^\star$. Conversely, if $\lambda \geq \tau^\star$, then, by $(i)$, the sequence $(t_n)$ is non increasing, greater than $\tau^\star$, and thus converges to $\ell \in [\tau^\star, \lambda]$ such that $\tau(\ell) = \ell$, and we also have $\ell = \tau^\star$. □

**Sto02 and FDR08.**

*Proof of Theorem 5.4.* As existence Condition C.4 holds, existence Condition C.11 also holds for procedure Sto02$(t)$, for any $t \in (0,1)$. Therefore, Theorem 4.15 ensures that the asymptotic threshold $\tau(t)$ of procedure Sto02$(t)$ is positive, and satisfies $G(\tau(t)) = \frac{\tau(t)}{\alpha}\overline{\pi_0}\,(t)$, where $\overline{\pi_0}\,(u) = \frac{1-G(u)}{1-u}$.

As uniqueness Condition C.3 and existence Conditions C.4 and C.7 hold, Theorem 4.7 ensures that the asymptotic threshold $\tau^\star$ of procedure FDR08 satisfies $\tau^\star \in (0,\kappa)$, and satisfies $G(\tau^\star) = f_\alpha(\tau^\star)$, where $f_\alpha : u \mapsto u/(\alpha + (1 -$ $\alpha)u)$ is the rejection curve of the FDR08 procedure. For any fixed $\lambda \in (t \wedge \kappa, 1)$, the FDR08($\lambda$) procedure defined by the capped threshold function

$$\mathcal{T}_\lambda(F) = \sup\{u \in [0, \lambda], F(u) \geq f_\alpha u\}$$

also has asymptotic threshold $\tau^\star$ according to Proposition 4.6, as $\lambda \geq \kappa$. For $F \in D[0,1]$ and $u \in [0, \lambda]$, let

$$c_\alpha(F, u) = \frac{u}{\alpha} \frac{1 - F(u)}{1 - u}.$$

As $G$ is concave, $u \mapsto \frac{1-G(u)}{1-u}$ is non increasing, so that $c_\alpha$ fulfills the requirements of Condition C.13. Therefore, as $F(u) \geq f_\alpha(u)$ may be written as $F(u) \geq c_\alpha(F, u)$, the result follows from the application of Theorem 5.2 to procedures FDR08($\lambda$) and Sto02($t$). □

**BKY06($\lambda$) and BR08($\lambda$).**

*Proof of Theorem 5.6.* As uniqueness Condition C.3 and existence Conditions C.8 and C.9 hold, Theorem 4.10 ensures that the asymptotic threshold $\tau^\star$ of procedure BR08($\lambda$) is the unique point in $(0, \lambda)$ such that $G(\tau^\star) = \frac{\tau^\star}{\alpha(1-\lambda)+\tau^\star}$, because the rejection curve $b_\alpha^\lambda$ of the BR08 procedure equals $\frac{u}{\alpha(1-\lambda)+u}$ for $u \leq \lambda$.

Existence Condition C.8 also ensures that $\tau(t)$ exists for any $t \leq \lambda$. For $F \in D[0,1]$ and $u \in [0, \lambda]$, let

$$c_\alpha(F, u) = \frac{u}{\alpha} \frac{1 - F(u)}{1 - \lambda}.$$

As $1 - G$ is non increasing, $c_\alpha$ fulfills the requirements of Condition C.13. Therefore, as $F(u) \geq b_\alpha(u)$ may be written as $F(u) \geq c_\alpha(F, u)$, the result follows from the application of Theorem 5.2 to procedures BR08($\lambda$) and BKY06($\lambda$). □

*Proof of Corollary 5.7.* According to the definition of $u_\lambda$ as the asymptotic threshold of the BH95 procedure at level $\lambda$, the asymptotic threshold $\tau^\star$ of procedure BR08($\lambda$) satisfies $\tau^\star \geq u_\lambda$ if and only if $G(\tau^\star) \leq \tau^\star$. According to the definition of the rejection curve $b_\alpha^\lambda$ of the BR08($\lambda$) procedure, this is equivalent to $\tau^\star/(\alpha(1-\lambda) + \tau^\star) \leq \tau^\star/\lambda$, that is, to $\tau^\star \geq \lambda - \alpha(1-\lambda)$. □

## Acknowledgments

I am extremely grateful to Stéphane Boucheron for his valuable help and useful advice throughout this work. I also wish to thank Étienne Roquain and Jean-Philippe Vert for insighful comments, and Antoine Chambaz for fruitful discussions. I am grateful to an Associate Editor for relevant suggestions that helped improve an earlier version of this manuscript.

## Notation used in the paper

TABLE 2
*Notation used throughout the paper*

| | |
|---|---|
| $\mathcal{H}_0$ | null hypothesis |
| $\mathcal{H}_1$ | alternative hypothesis |
| $m$ | number of tested hypotheses |
| $m_0(m)$ | number of tested hypotheses |
| $\pi_0 = m_0(m)/m$ | proportion of true null hypotheses |
| $(P_i)_{i\in\{1\dots m\}}$ | associated $p$-values |
| $\mathcal{M}$ | Multiple Testing Procedure (MTP), cf. Definition 2.1 |
| $R_m$ | number of rejections among $m$ hypotheses for a given MTP |
| $V_m$ | number of false rejections among $m$ hypotheses for a given MTP |
| $\mathsf{FDP} = \frac{V_m}{R_m \vee 1}$ | FDP attained by a given MTP |
| $\mathsf{FDR} = \mathbb{E}\left[\frac{V_m}{R_m \vee 1}\right]$ | FDR attained by a given MTP |
| $\mathcal{T}$ | threshold function of a MTP, cf. Definition 2.2 |
| $r_\alpha = r(\alpha, \cdot)$ | rejection curve of a MTP |
| $\mathcal{A}$ | level function of a MTP |
| $\mathcal{U}$ | $(F,\alpha) \mapsto \sup\{u \in [0,1], F(u) \geq r_\alpha(u)\}$ |
| $\mathsf{p}(t) = \mathsf{pFDR}(t)$ | positive pFDR attained at $t$ |
| $\dot{\mathsf{p}}(t) = \frac{d}{dt}(\mathsf{pFDR})(t)$ | |
| $G_0 : u \mapsto u\mathbf{1}_{[0,1]}(u)$ | distribution function of $p$-values under $\mathcal{H}_0$ (Uniform) |
| $G_1$ | distribution function of $p$-values under $\mathcal{H}_1$ ($C^1$, concave) |
| $g_1 = G_1'$ | density of $p$-values under the alternative hypothesis |
| $G = \pi_0 G_0 + (1-\pi_0)G_1$ | |
| $g = \pi_0 + (1-\pi_0)g_1$ | |
| $\widehat{\mathbb{G}}_{0,m}(t) = \frac{\sum_{\{i/\mathcal{H}_0 \text{ true}\}} \mathbf{1}_{P_i \leq t}}{m_0(m)}$ | empirical distribution function of $p$-values under $\mathcal{H}_0$ (unobservable) |
| $\widehat{\mathbb{G}}_{1,m}(t) = \frac{\sum_{\{i/\mathcal{H}_1 \text{ true}\}} \mathbf{1}_{P_i \leq t}}{m-m_0(m)}$ | empirical distribution function of $p$-values under $\mathcal{H}_1$ (unobservable) |
| $\widehat{\mathbb{G}}_m = \pi_0\widehat{\mathbb{G}}_{0,m} + \widehat{\mathbb{G}}_{1,m}$ | empirical distribution function of $p$-values (observable) |
| $\mathbb{B}$ | standard Brownian bridge on $[0,1]$ |
| $\mathbb{Z}_0 \overset{(d)}{=} \mathbb{B}$ | limit in distribution of $\sqrt{m}(\widehat{\mathbb{G}}_{0,m} - G_0)$ |
| $\mathbb{Z}_1 \overset{(d)}{=} \mathbb{B} \circ G_1$ | limit in distribution of $\sqrt{m}(\widehat{\mathbb{G}}_{1,m} - G_1)$ |
| $\mathbb{Z} = \pi_0\mathbb{Z}_0 + (1-\pi_0)\mathbb{Z}_1$ | limit in distribution of $\sqrt{m}(\widehat{\mathbb{G}}_m - G)$ |
| $\tau^\star = \mathcal{T}(G)$ | asymptotic threshold of the MTP with threshold function $\mathcal{T}$ |
| $\alpha^\star = \lim_{u\to 0} \frac{1}{g(u)}$ | critical value of the BH95 procedure, as defined by [5] |
| $A \gg B$ | procedure $A$ is more powerful than procedure $B$ |